\documentclass{article}
\usepackage{arxiv}
\usepackage[utf8]{inputenc} 
\usepackage[T1]{fontenc}    
\usepackage{nicefrac}       
\usepackage{microtype} 
\usepackage{multirow}
\usepackage{amsmath,amssymb,psfrag}
\usepackage{graphicx}
\usepackage{fullpage}
\usepackage{listings}
\usepackage{paralist}
\usepackage{appendix}
\usepackage{amsfonts}
\usepackage{comment} 
\usepackage{booktabs} 
\usepackage{cancel}
\usepackage{helvet}
\usepackage{tikz,tikz-cd}
\usepackage{pgfplots,pgfplotstable}
\usepackage{lineno}
\usepackage{hyperref}
\usepackage{mathtools}
\usepackage{amsthm}
\usepackage{enumerate}
\usepackage{enumitem}
\usepackage{subfig,color}
\usepackage{xcolor}
\usepackage{epstopdf}
\usepackage{bm}
\usetikzlibrary{pgfplots.groupplots}
\usetikzlibrary{positioning}
\usetikzlibrary[shapes,arrows,trees]
\usetikzlibrary{matrix,decorations.pathmorphing}

\newcommand{\xddots}{%
  \raise 4pt \hbox {.}
  \mkern 6mu
  \raise 1pt \hbox {.}
  \mkern 6mu
  \raise -2pt \hbox {.}
}

\newcommand{\bx}{\bm{x}}

\newtheorem{theorem}{Theorem}
\newtheorem{corollary}{Corollary}[theorem]
\newtheorem{lemma}{Lemma}[theorem]

\newtheoremstyle{sltheorem}
{}                
{}                
{\slshape}        
{}                
{\bfseries}       
{.}               
{ }               
{}                
\theoremstyle{sltheorem}

\newtheorem{theorem2}{Theorem}
\newtheorem{definition}{Definition}
\numberwithin{equation}{section}
\numberwithin{figure}{section}
\numberwithin{table}{section}

\numberwithin{equation}{section}
\numberwithin{figure}{section}
\numberwithin{table}{section}
\newcommand{\myref}[1]{~(\ref{#1})}

\title{OS-net: Orbitally Stable Neural Networks}

\author{
  Marieme Ngom \\
  Mathematics and Computer Science division\\
  Argonne National Laboratory\\
  Lemont, IL 60439 \\
  \texttt{mngom@anl.gov} \\
   \And
Carlo Graziani \\
  Mathematics and Computer Science division\\
  Argonne National Laboratory\\
  Lemont, IL 60439 \\
  \texttt{cgraziani@anl.gov} \\
}
\pgfplotsset{compat=1.18} 
\begin{document}

\maketitle

\begin{abstract}
We introduce OS-net (Orbitally Stable neural NETworks), a new family of neural network architectures specifically designed for periodic dynamical data. OS-net is a special case of Neural Ordinary Differential Equations (NODEs) and takes full advantage of the adjoint method based backpropagation method. Utilizing ODE theory, we derive conditions on the network weights to ensure stability of the resulting dynamics. We demonstrate the efficacy of our approach by applying OS-net to discover the dynamics underlying the R\"{o}ssler and Sprott's systems, two dynamical systems known for their period doubling attractors and chaotic behavior. 
\end{abstract}
\section{Introduction}
The study of periodic orbits of systems of the form
\begin{equation}\label{eq:eq_diff_gen}
    \dot{\bm x} = f(\bm x),\; \bm x(0) = \bm x_0, \; \bm x \in \mathbf{U} \subset \mathbf{R}^n 
\end{equation}
is an important area of research within the field of nonlinear dynamics with applications in both the physical (astronomy, meteorology) and the nonphysical (economics, social psychology) sciences. In particular, periodic orbits play a significant role in chaos theory. In \cite{Devaney2003}, chaotic systems are defined as systems that are sensitive to initial conditions, are topologically transitive (meaning that any region of the phase space can be reached from any other region), and have dense periodic orbits. Notably, chaotic systems are constituted of infinitely many Unstable Periodic Orbits (UPOs) which essentially form a structured framework, or a "skeleton", for chaotic attractors. A periodic orbit is (orbitally) unstable if trajectories that start near the orbit do not remain close to it. Finding and stabilizing UPOs is an interesting and relevant research field with numerous applications such as the design of lasers \cite{Roy1992}, the control of seizure activities \cite{Schiff1994} or the design of control systems for satellites \cite{Wiesel1993}. An important tool when studying the stability of periodic orbits of a given system is the Poincar\'e or return map which allows one to study the dynamics of this system in a lower dimensional subspace. It is well-known that the stability of a periodic orbit containing a point $\bm x_0$ is inherently connected to the stability of $\bm x_0$ as a fixed point of the corresponding Poincar\'e map. However, explicitly computing Poincaré maps has been proven to be highly challenging and inefficient \cite{Teschl2012}. With the emergence of data-driven approaches, researchers in \cite{Bramburger2020} proposed a data-driven computation of Poincaré maps using the SINDy method \cite{Brunton2016}. Subsequently, they leveraged this technique in to develop a method for stabilizing UPOs of chaotic systems  \cite{Bramburger2021}. \\
As a matter of fact, researchers have been increasingly exploring the intersection of machine learning and differential equations in recent years. For example, Partial Differential Equations (PDEs) inspired neural network architectures have been developed for image classification in \cite{Ruthotto2019, Sun2020}. On the other hand, data-driven-based PDE solvers were proposed in \cite{Sirignano2018} while machine learning has been effectively utilized to discover hidden dynamics from data in \cite{raissi2019Pinn, Brunton2016, Schaeffer2017}. One notable example of such intersectional work is Neural Ordinary Differential Equations (NODEs), which were introduced in \cite{Chen2018}. NODEs are equivalent to continuous residual networks that can be viewed as discretized ODEs \cite{Haber2018}. This innovative approach has led to several extensions that leverage well-established ODE theory and methods \cite{Dupont2019, Zhang2022, Ott2021, zhuang2020, Yan2020, Haber2018} to develop more stable, computationally efficient, and generalizable architectures.\\
In the present work, we aim at learning dynamics obtained from chaotic systems with a shallow network featuring a single hidden layer, wherein the network's output serves as a solution to the dynamical system
\begin{equation}\label{eq:auton}
    \Dot{\bx} = \bm W_d^T\sigma(\bm W_e^T \bx + \bm b_e) , \;\; \bx(0) = \bx_0.
\end{equation}
 where $\bm W_e$ are the input-to-hidden layer weights, $\bm b_e$ the corresponding bias term, $\bm W_d$ the hidden-to-output layer weights, and $\sigma$ the activation function of the hidden layer. The proposed network is a specific case of NODEs and fully utilizes the adjoint method-based \cite{Pontryagin1987} weight update strategy introduced in \cite{Chen2018}. Our primary objective is to establish sufficient conditions on the network parameters to ensure that the resulting dynamics are orbitally stable. We base our argument on the finding that the stability of Poincaré maps is equivalent to the stability of the first variational equation associated with the dynamical system under consideration \cite{Teschl2012}. We then build on the stability results of linear canonical systems presented in \cite{Krein1983} to derive a new regularization strategy that depends on the matrix $\bm J = \bm W_e^T \bm W_d^T$ and not on the weight matrices taken independently. We name the constructed network OS-net for Orbitally Stable neural NETworks.\\
Since we are dealing with periodic data, the choice of activation function is critical. Indeed, popular activation functions such as the $\mathrm{sigmoid}$ or the $\tanh$ functions do not preserve periodicity outside the training region. A natural choice would be sinusoidal activations however these do not hold desired properties such as monotonicity. Furthermore, they perform poorly \cite{Parascandolo2017} on the training phase because the optimization can stagnate in a local minimum because of the oscillating nature of sinusoidal functions. In \cite{Ngom2021}, the authors constructed a Fourier neural network (i.e a neural network that mimics the Fourier Decomposition) \cite{Silvescu1999, Zhumekonov2019} that uses a $\sin$ activation but had to enforce the periodicity in the loss function to ensure that periodicity is conserved outside of the training region. The activation functions $x + \frac{1}{a}\sin^2 (ax)$ -called snake function with frequency $a$- and $x + \sin x$ were proposed in \cite{Ziyin2020} for periodic data and were proven to be particularly well suited to periodic data. As such, we use both these activation functions in this work. \\
This paper is organized as follows: in section\myref{sec:bkgd} we present the OS-net's architecture and the accompanying new regularization strategy. In section\myref{sec:results} we showcase its performance on simulated data from the chaotic R\"{o}ssler and Sprott systems and perform an ablation study to assess the contributions of the different parts of OS-net.
\section{Building OS-net}\label{sec:bkgd}
\subsection{Background}
In this chapter, we recall the main results on the stability of periodic orbits of dynamical systems we will be using to build OS-net. We refer readers to the appendices for more details about orbits of dynamical systems.\\
We consider the system 
\begin{equation}\label{eq:simple_dyn}
    \dot{\bm x} = f(\bm x),\;\; \bm x(0) = \bm x_0.
\end{equation}
and suppose it has a periodic solution $\phi(t, x_0)$ of period $T$. We denote $\gamma(x_0)$ a   periodic orbit corresponding to $\phi(t, x_0)$.  Stability of periodic orbits have been widely studied in the literature. It is, in particular well-known \cite[Chapter~12]{Teschl2012} that the stability of periodic orbits of Equation\myref{eq:simple_dyn} is linked to the stabiity of its First Variational (FV) problem
\begin{equation}\label{eq:fvar}
    \Dot{\bm y} = \bm A(t)\bm y, \;\;\bm A(t) = d\left(f(\bm x)\right)_{(\Phi(t,\bm x_0))}\;\; \text{and } \bm A(t+T) = \bm A(t).
\end{equation}
which is obtained by taking the gradient of Equation\myref{eq:simple_dyn} with respect to $x$ at $\phi(t, x_0)$. As such, the first variational problem describres the dynamics of the state variable $y = d \left(\phi(t, t_0, \bm x)\right)$ and is a linear system as the matrix $\bm A(t)$ does not depend on $y$. \\
To assess the stability of OS-net, we investigate the first variational equation associated with Equation\myref{eq:auton}. OS-net's FV is given by
\begin{equation*}
    \Dot{\bm y} = \bm W_d^T diag\left(\sigma^{'}\left(\bm W_e^T \phi(t,\bm x_0) + \bm b_e\right)\right) \bm W_e^T\bm y,
\end{equation*}
and if we make the change of variables $\bm z = W_e^T\bm y$, this equation becomes
\begin{equation}\label{eq:osnet_fvar}
    \Dot{\bm z} = \bm J \bm H(t) \bm z,
\end{equation}
where $J = \bm W_e^T \bm W_d^T$ and $\bm H(t)  = diag\left(\sigma^{'}\left(\bm W_e^T \phi(t,\bm x_0) + \bm b_e\right)\right)$ is periodic. This formulation can be seen as a generalization of linear canonical systems with periodic coefficients 
\begin{equation}\label{eq:krein_cano}
    \Dot{\bm y} = \lambda \bm J_m \bm H(t)\bm y
\end{equation}
where $$\bm J_m = \begin{pmatrix}
0 & \bm I_m\\
-\bm I_m & 0
\end{pmatrix}, \bm I_m \text{ is the identity matrix of size $m$}, $$
$\bm H$ is a periodic matrix-valued function and $\lambda \in \mathrm{R}$. Stability of such systems was extensively studied in \cite{KreinJakub1983} and in particular in \cite{Krein1983}. We recall the main definitions and results from \cite{Krein1983} and build upon these to  derive stability conditions for OS-net. In particular, we give the definition of stability zones for Equation\myref{eq:krein_cano} and provide the main stability results we will base our study on.
\begin{definition}
A point $\lambda = \lambda_0\;(-\infty < \lambda_0 < \infty)$ is called a $\lambda$-point of stability of Equation \myref{eq:krein_cano} if for $\lambda = \lambda_0$ all solutions of Equation \myref{eq:krein_cano} are bounded on the entire t-axis.\\
If, in addition, for $\lambda = \lambda_0$ all solutions for any equation 
\begin{equation*}
    \Dot{\bm y} = \lambda \bm J_m \bm H_1(t)\bm y
\end{equation*}
with a periodic symmetric. matrix valued function $H_1(t) = H_1(t+T)$ sufficiently close to $H(t)$ are bounded, then $\lambda_0$ is a $\lambda$-point of strong stability of Equation \ref{eq:krein_cano}.
\end{definition}
The set of $\lambda$-point of strong stability of Equation \ref{eq:krein_cano} is an open-set that decomposes into a system of disjoint open intervals called $\lambda$-zones of stability of Equation\ref{eq:krein_cano}. If a zone of stability contains the point $\lambda = 0$ then it is called a central zone of stability.
\begin{definition}
    We say that Equation \myref{eq:krein_cano} is of positive type if 
    $$
    \bm H \in \mathrm{P}_n(T) = \{ \bm A(t) \; symmetric \; s.t \; \bm A(t) \geq 0 \; (0 \leq t \leq T) \; \text{and} \int_0^T \bm A(t)dt >0 \}.
    $$
 $\bm A(t) \geq 0$ means $\forall x \in \mathbf{R}^n, \langle \bm A(t) \bm x,\; \bm x \rangle \geq 0$ and $\int_0^T \bm A(t)dt >0$ means $ \int_0^T \langle \bm H(t) \bm x,\; \bm x \rangle dt >0$
\end{definition}
\begin{definition}
    Let $\bm A$ be a square matrix with non-negative elements. We denote by $\mathcal{M}(\bm A)$ the least positive eigenvalue among its eigenvalues of largest modulus. Note that Perron's theorem (1907) guarantees the existence of $\mathcal{M}(\bm A)$ \cite{horn2012matrix}.
\end{definition}
We can now state the main result we will derive our regularization from:
\begin{theorem}[\cite{Krein1983} section 7, criterion $I_n$]\label{theo:krein}
    A real $\lambda$ belongs to the central zone of stability of an Equation \myref{eq:krein_cano} of positive type, if 
$$|\lambda| < 2 \mathcal{M}^{-1} (\bm C)$$
where $\bm C = \bm J_{m_{a}}\int_0^T \bm H_a(t)$. 
If $K$ is a matrix, $K_a$ is the matrix obtained by replacing the elements of $K$ by their absolute values.
\end{theorem}
The proof of this theorem is recalled in Appendix\myref{appendixB}.

\subsection{Architecture and stability of OS-net }\label{sec:os-net}
To base the stability of OS-net on stability theory for systems of type Equation\myref{eq:krein_cano}, we need the matrix-valued function $\bm H(t)$ and the matrix $\bm J$ in Equation \myref{eq:osnet_fvar} to be respectively of positive type and skew-symmetric.
 \begin{figure}[htb]
 \vspace{-20mm}
    \centering
    \subfloat{\includegraphics[width=0.45\textwidth]{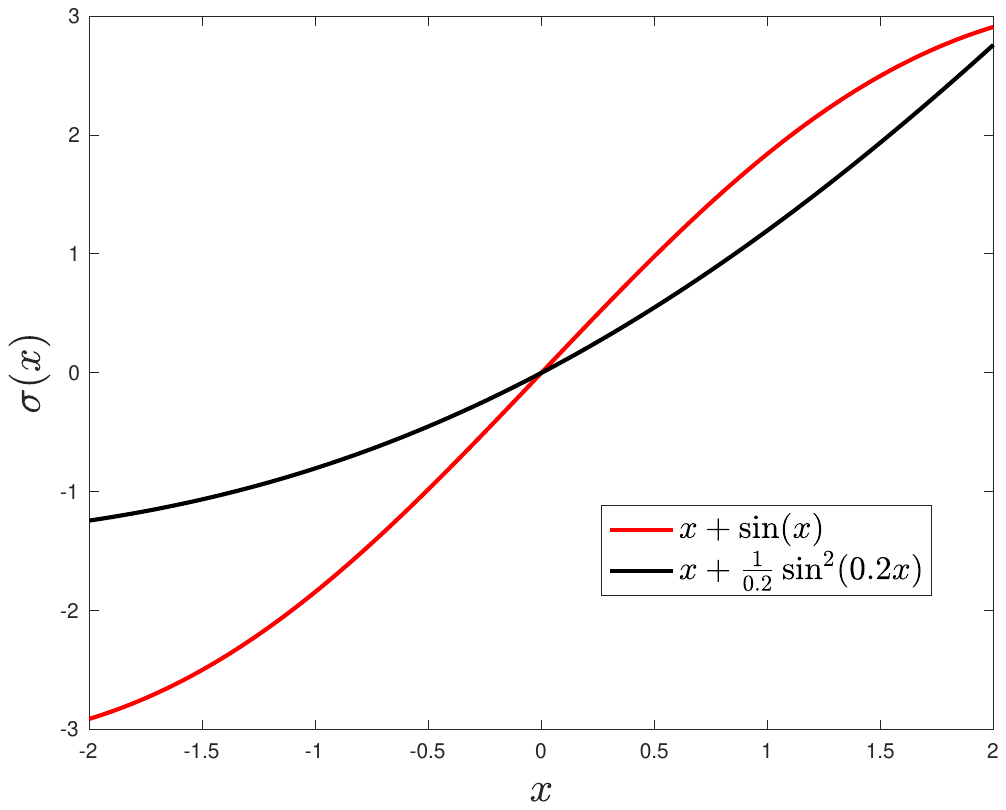}}
     \subfloat{\includegraphics[width=0.45\textwidth]{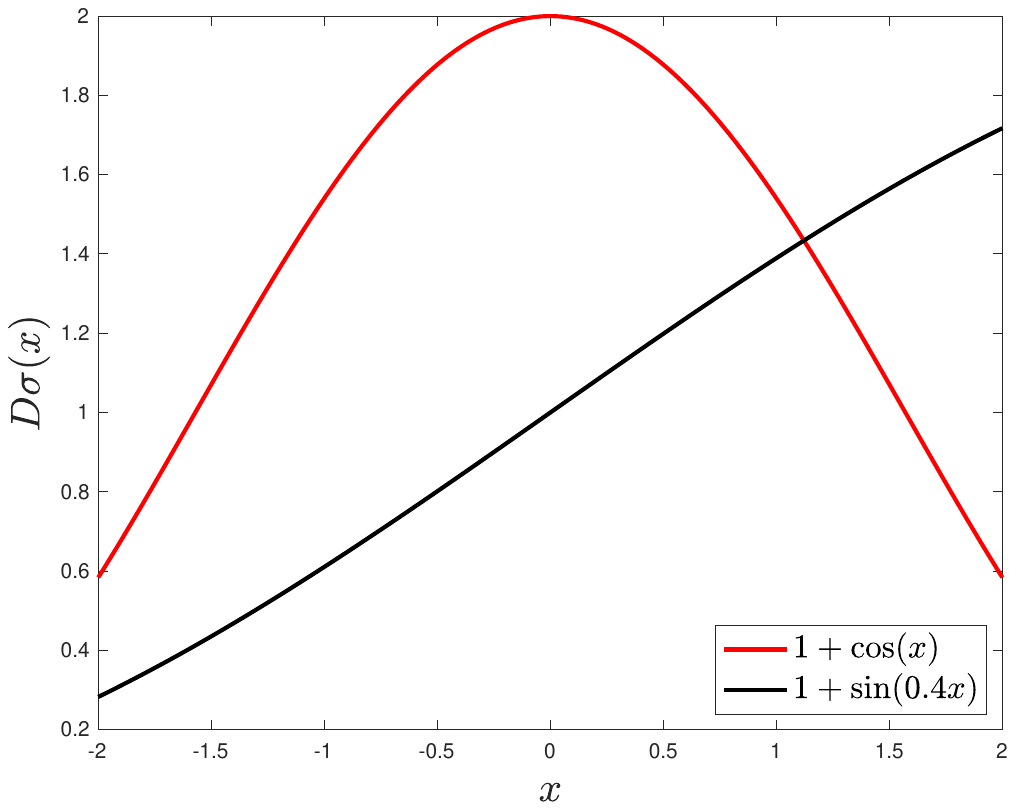}}
     \vspace{-20mm}
     \caption{\;Left: Activation functions $x + \sin(x)$ and $x + \frac{1}{0.2} \sin^2(0.2 x)$. Right: Derivatives of the activation functions}\label{fig:activ}
\end{figure}
To ensure $\bm H(t)$ is of positive type, it is sufficient to use activation functions that are increasing since they have positive derivatives and diagonal matrices with positive elements are of positive type. Fortunately, many common activation functions ($\tanh$ or $sigmoid$) have that property. In this paper,  we use the  strictly monotonic activation functions $x + \sin (x)$ and (the snake function) $x + \frac{1}{a} \sin^2(a x), \; a \in \mathbf{R}$ displayed in Figure\myref{fig:activ}. These activation functions were proved to be able to learn and extrapolate periodic functions in \cite{Ziyin2020}.

Let us now pay attention to the matrix $\bm J = \bm W_e^T \bm W_d^T$. To ensure $J$ is skew-symmetric, we introduce the matrices $\bm W \in \mathrm{M}_{n, 2m}(\mathbf{R})$ and $\bm K \in \mathrm{M}_{2m}(\mathbf{R})$ where $n$ is the input size (i.e the size of $\bm x$) and $2m$ the size of the hidden layer. We then set $\bm W_e^T = \bm W$, $\bm W_d^T = \bm \Omega \bm W^T$, and $\bm \Omega = \bm K - \bm K^T$. Note that the size of the hidden layer which is the size of $\Omega$ needs to be even. Otherwise, $\bm \Omega$ would be a singular matrix. The elements of the matrices $\bm W$ and $\bm K$ are the hyperparameters of the network that will be optimized during training. Now, knowing that any real skew-symmetric matrix $\bm J$ is congruent to $\bm J_m$ \cite{Yakubovich1975}, there exists a real invertible matrix $\bm S$ such that $$\bm J_m = \bm S^T \bm J \bm S$$ and Equation\myref{eq:osnet_fvar} is equivalent to Equation\myref{eq:krein_cano}. In fact, let $\bm u = \bm S^T \bm z$ in Equation\myref{eq:osnet_fvar}, we obtain
$$\dot{\bm u} = \lambda \bm J_m \bm S^{-1} \bm H(t) (\bm S^{-1})^T \bm u = \lambda \bm J_m \Tilde{\bm H}(t) \bm u$$
where $\Tilde{\bm H}(t) = \bm S^{-1} \bm H(t)  (\bm S^{-1})^T \in \mathrm{P}_n(T)$. We can now apply Theorem\myref{theo:krein} to OS-net and state that OS-net is stable if 
\begin{equation}
    1 < 2 \mathcal{M}^{-1}\left(\bm J_a \int_{0}^{T} \bm H(t)dt\right).
\end{equation}
Note that since $\bm H(t)$ is a diagonal matrix with positive elements, $\bm H_a(t) = \bm H(t)$.
We can now prove the following result that will justify our regularization strategy:
\begin{corollary}
Suppose the activation function $\sigma$ is strictly increasing with a uniformly bounded derivative. Then, OS-net is stable if 
\begin{equation}
    ||\bm J_a||_2 < \frac{2}{L T} 
\end{equation}\label{eq:osnet_stab_cond}
where $L$ is the superior bound of the derivative of the activation function.
\end{corollary}
\begin{proof}
Let $\mu$ be any eigenvalue of $\bm J_a \int_{0}^{T} \bm H(t)dt$ then $|\mu| \leq  \left\Vert \bm J_a \int_{0}^{T} \bm H(t)dt \right\Vert_2$. Knowing that any norm in $\mathbf{R}^{n,n}$ can be rescaled to be submultiplicative (i.e. $||AB||_2 \leq ||A||_2||B||_2$), we obtain $$|\mu| \leq  \Vert \bm J_a\Vert_2 \left\Vert\int_{0}^{T} \bm H(t)dt  \right\Vert_2$$
which leads to 
$$\mathcal{M}\left( \bm J_a \int_{0}^{T} \bm H(t)dt \right) \leq || \bm J_a||_2 \left\Vert\int_{0}^{T} \bm H(t)dt  \right\Vert_2$$
If $L$ is the superior bound of the derivative of the activation function, then, since $\bm H(t)$ is a diagonal matrix, we have $\mathcal{M}\left( \bm J_a \int_{0}^{T} \bm H(t)dt \right) \leq L T ||\bm J_a||_2 $. Therefore, OS-net is stable if $1 < \frac{2}{L T} ||\bm J_a||_2^{-1}$
i.e $||\bm J_a||_2 < \frac{2}{L T}.$
\end{proof}

All in all our minimization problem becomes
\begin{equation}\label{eq:loss_pre_final}
  L(\bm x_o, g(f(\bm x_0))) = ||g(f(\bm x_0)) - \bm x_0||^2_2\;\; \text{s.t.}\;\; ||\bm J_a||_2 < \frac{2}{L T}
\end{equation}
and this formulation is equivalent \cite{Bach2011} to 
\begin{equation}\label{eq:loss_final}
  L(\bm x_0, g(f(\bm x_0))) = ||g(f(\bm x_0)) - \bm x_0||^2_2 + \alpha ||\bm J_a||^2_2
\end{equation}
where $\alpha \in \mathbf{R}$ can be fine-tuned using cross-validation. We thus have derived a new regularization strategy that stabilizes the network. By controlling the norm of $J_a = \vert\bm W_e^T \bm W_d^T \vert $, we ensure solutions of Equation\myref{eq:osnet_fvar} and consequently periodic orbits of Equation\myref{eq:auton} are stable. We validate these claims in the next section by running a battery of tests on simulated data from dynamical systems known for their chaotic behavior.


\section{Numerical results}\label{sec:results}
In this section, we showcace the learning capabilities and stability of OS-net on different regimes of the R\"{o}ssler \cite{Rossler1976397} and of the Sprott systems \cite{Sprott1997}. In all of the following experiments, the data was generated using Matlab's ode45 solver. We take snapshots at different time intervals to obtain the  data used to train OS-net.
.

We used the $LBFGS$ optimizer with a learning rate $lr = 1.$ and the strong Wolfe \cite{Wolfe1969} line search algorithm for all the experiments. Our code uses Pytorch and
all the tests were performed on a single GPU\footnote[1]{We base our code on the Neural ode implementation in \cite{msurtsukov_github} } using Argonne Leadership Computing Facility (ALCF)'s Theta/ThetaGPU \cite{ThetaGPU}.

\subsection{The R\"{o}ssler system}\label{sec:ross}
As in \cite{Bramburger2021}, we consider the R\"{o}ssler system
\begin{align}
    \dot{x} &= -y - z\\
    \dot{y} &= x + 0.1 y \\
    \dot{z} &= 0.1 + z (x - c)
\end{align}\label{eq:Rossler}
where $c \in \mathbf{R}$. R\"{o}ssler introduced this system as an example of simple chaotic system with a single nonlinear term ($zx$). As $c$ increases, this system displays period doubling bifurcations leading to chaotic behavior. Here, we consider the values $c = 6$ and $c = 18$. 
\subsubsection{c = 6, period-2 attractor}\label{sec:ross6}
First, we set $c = 6$ and initialize the trajectory at $[ x_0, y_0, z_0] = [0, -9.1238 ,0]$. In this regime, the R\"{o}ssler system possesses a period-2 attractor \cite{Bramburger2021}. 
 \begin{figure}[htb]
    \centering   \subfloat{\includegraphics[width=0.45\textwidth]{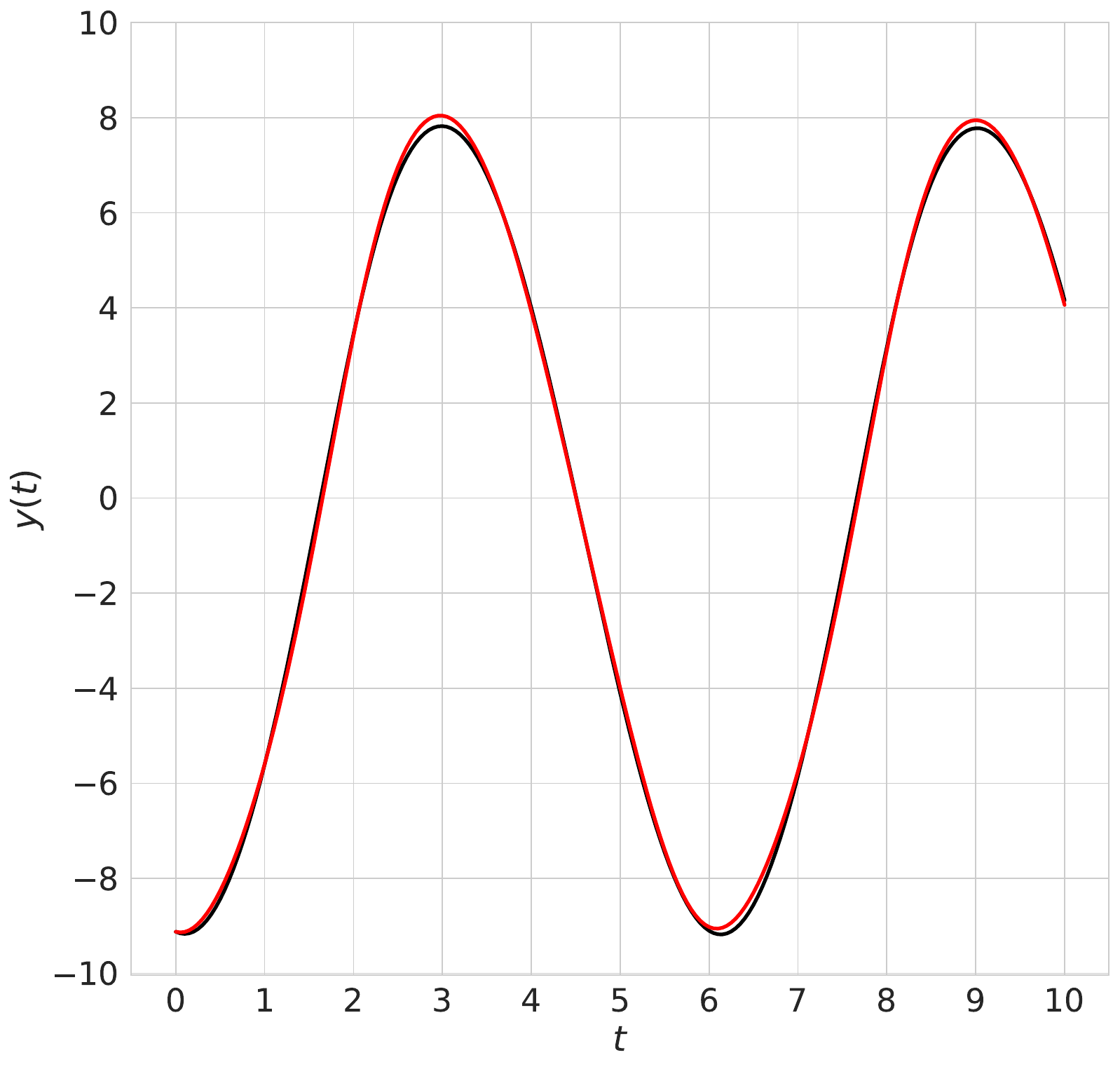}}
    \hspace{1mm}
     \subfloat{\includegraphics[width=0.45\textwidth]{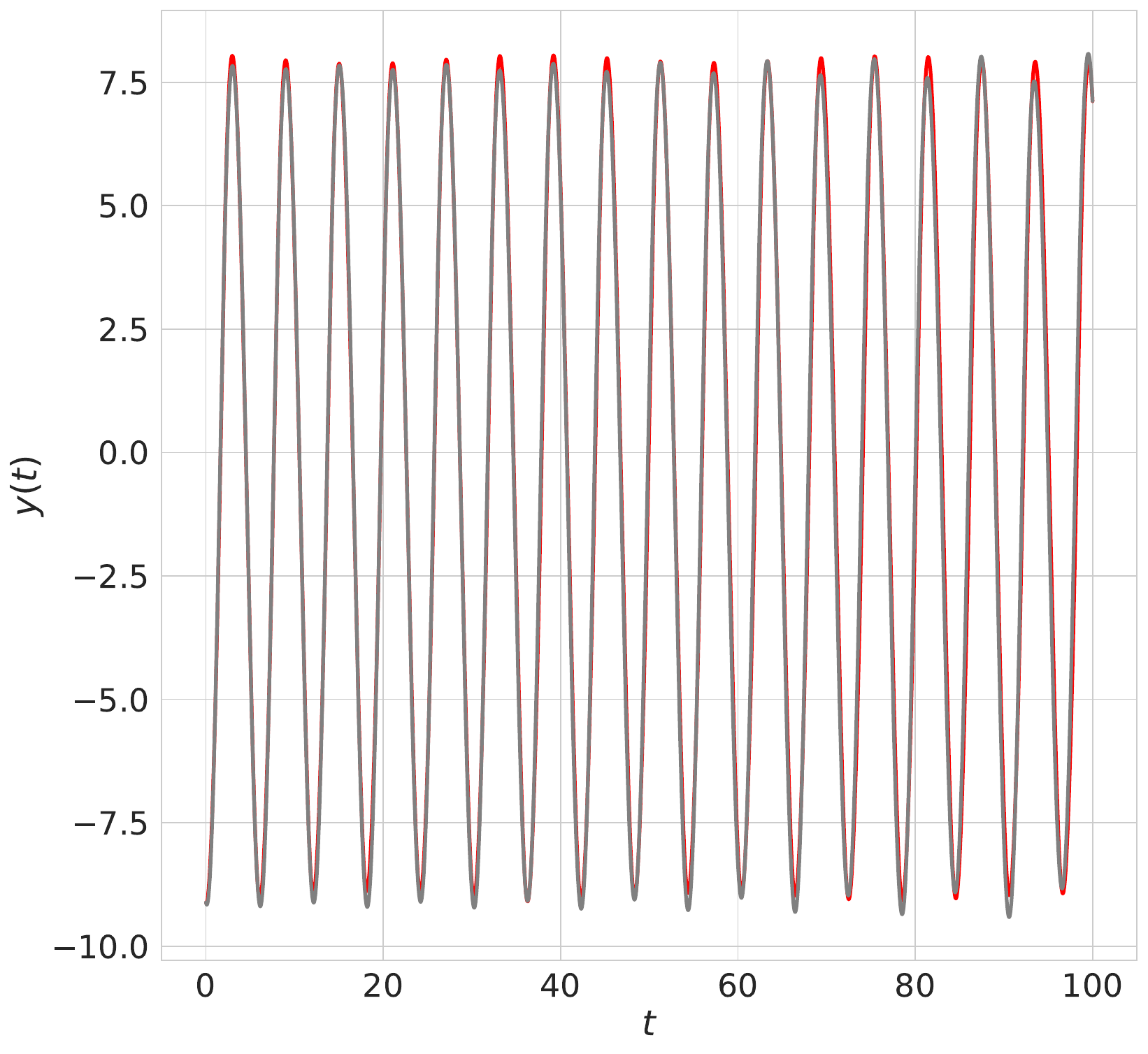}}
     \caption{\;Left: Training data in black along with the learned dynamics in red in the training time interval $[0, 10]$. Right: Target dynamics in gray along with the data generated by OS-net on the time interval $[0, 100].$ }\label{fig:ross_train_pred}
\end{figure}
\begin{table}[htb]
  \caption{Norm of $J_a$}
  \label{tab:normJ}
  \centering
  \begin{tabular}{lll}
    \toprule
    \multicolumn{2}{c}{Part}                   \\
    \cmidrule(r){1-2}
    System     & Attractor type     & $||J_a||$ \\
    \midrule
    R\"{o}ssler, $c = 6$ & Period-2  & 0.9937     \\
    R\"{o}ssler, $c = 18$    & Chaotic & 0.6318     \\
    Sprott, $\mu = 2.1$     & Period-2       & 0.0085 \\
    \bottomrule
  \end{tabular}
\end{table}

 \begin{figure}[htb]
    \centering
    \subfloat{\includegraphics[width=0.45\textwidth]{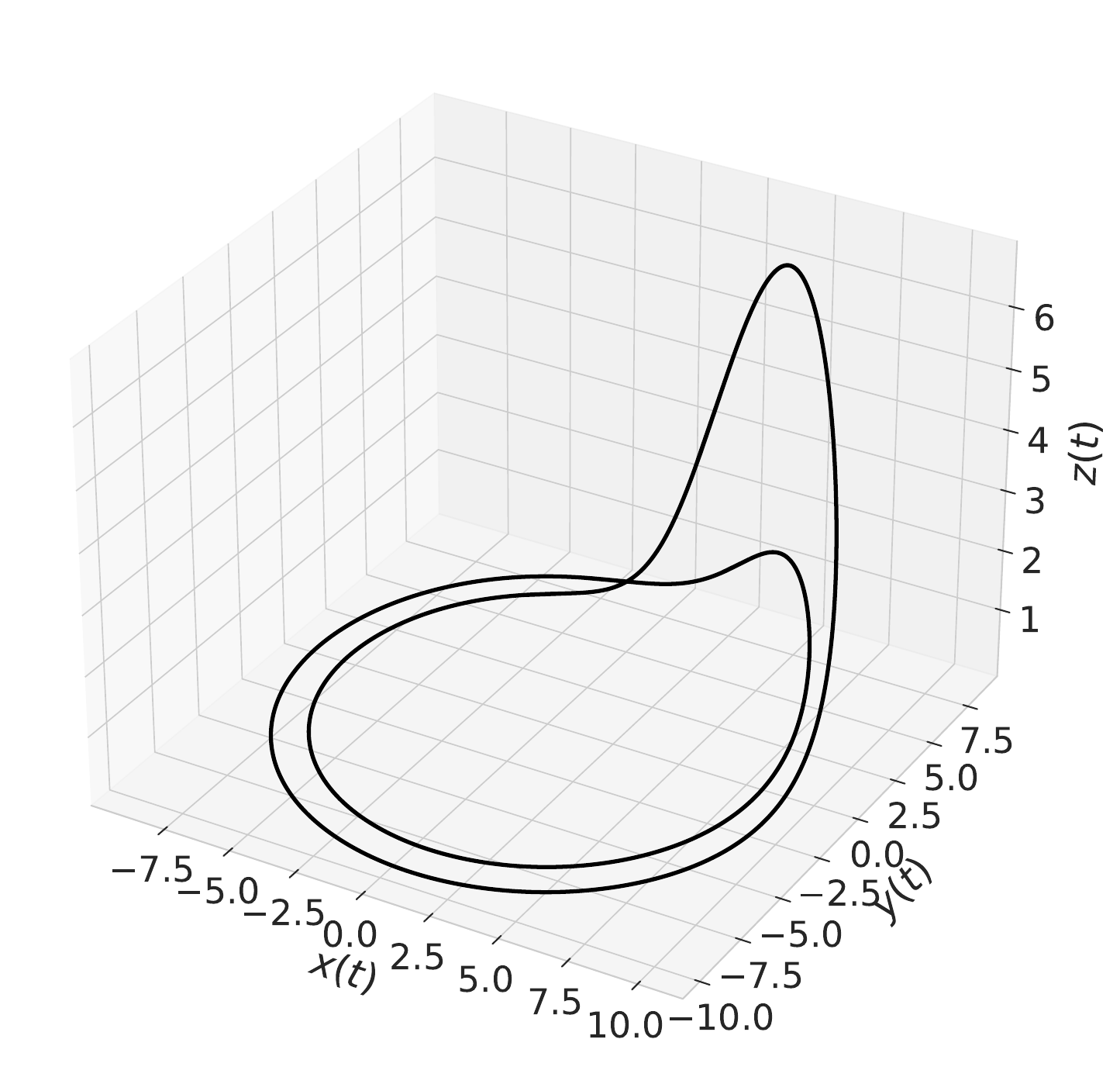}}
     \subfloat{\includegraphics[width=0.45\textwidth]{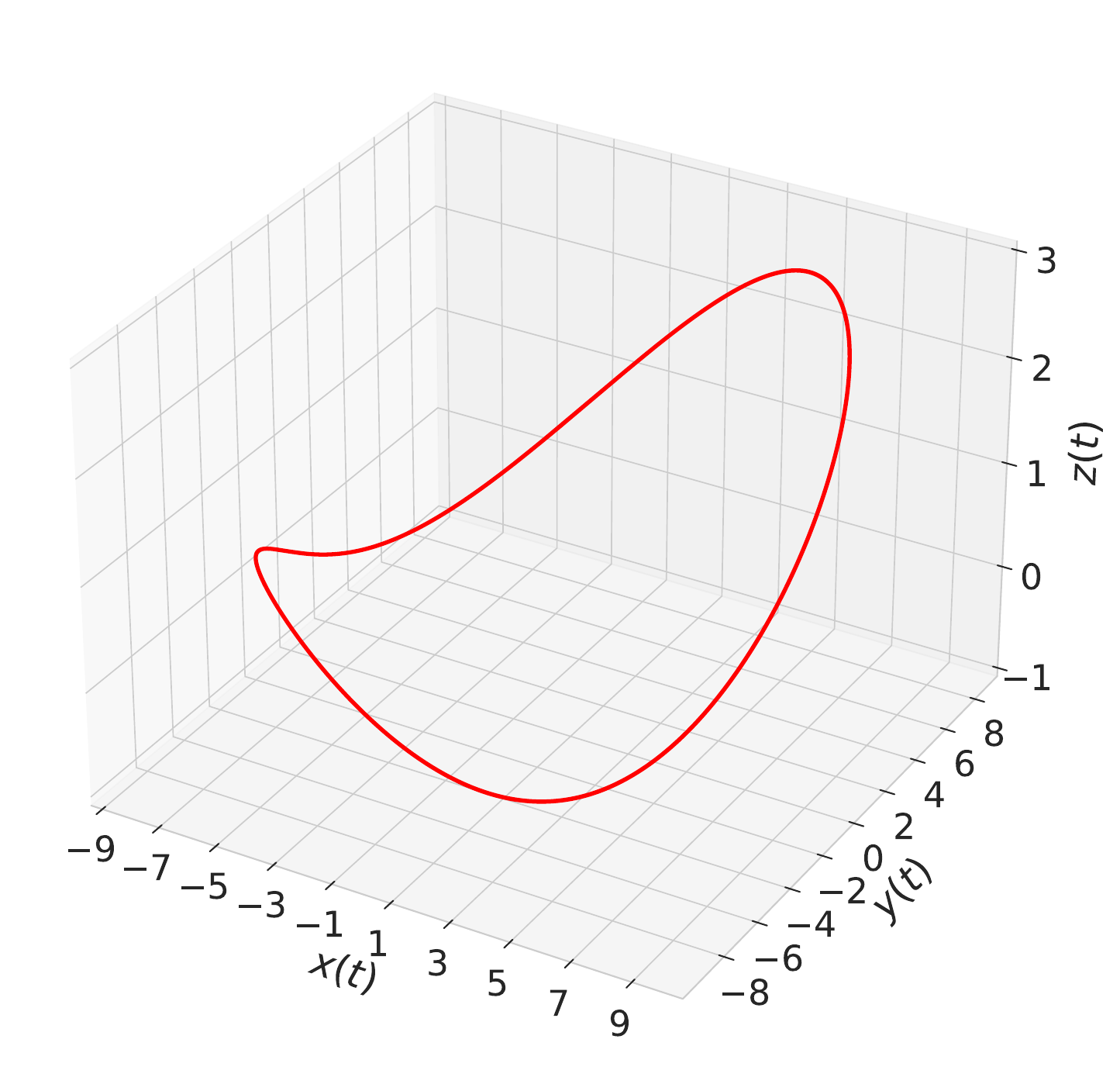}}
     \caption{\;Left: R\"{o}ssler's period-2 attractor. Right: stable OS-net period-1  attractor (right)}\label{fig:ross_osnet}
\end{figure}
\begin{figure}[htb]
     \centering
     \includegraphics[width=0.45\textwidth]{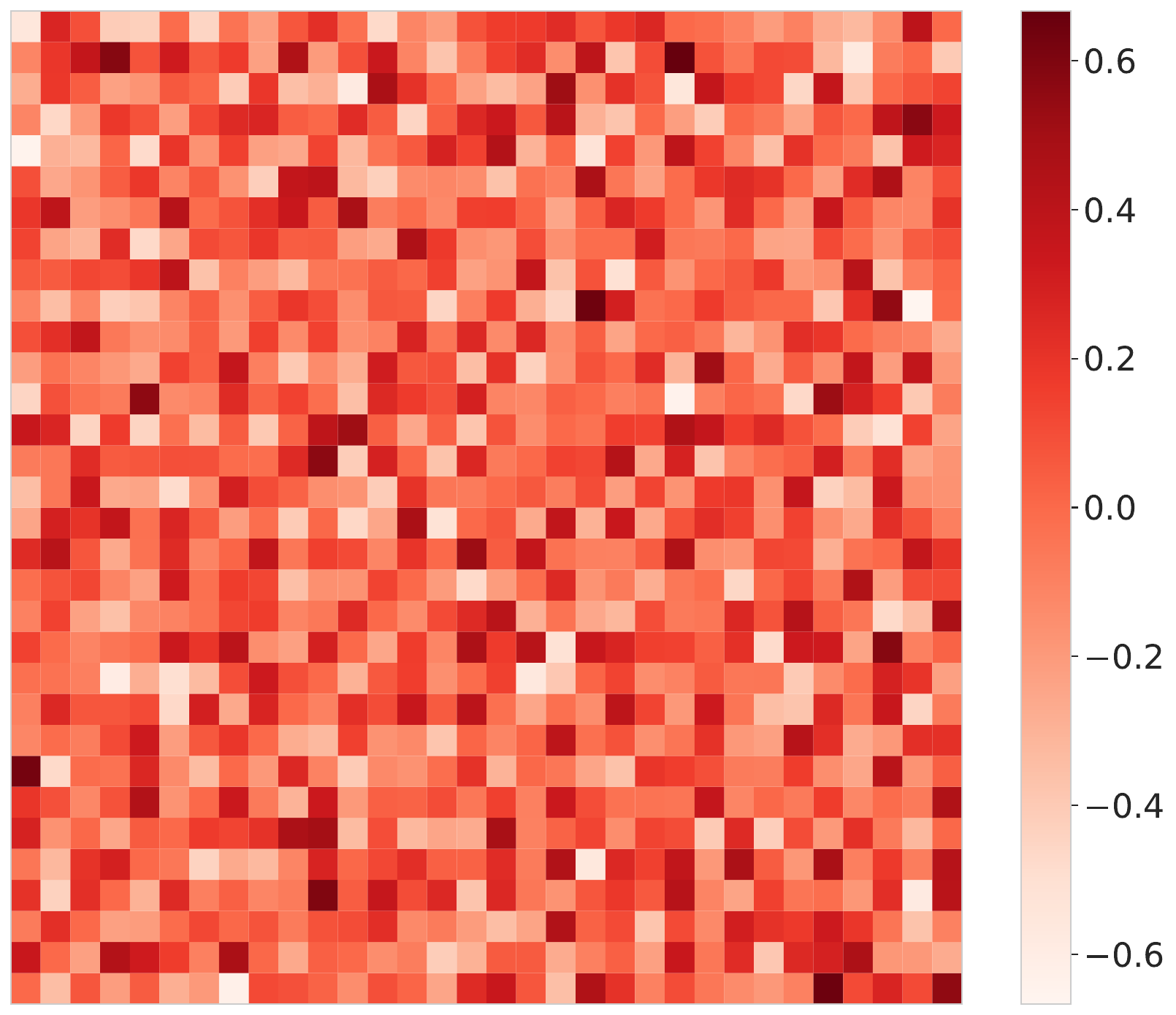}
\caption{\;Elements of the matrix $\Omega$}\label{fig:colormap_ross}
 \end{figure}
\begin{figure}[htb]
    \centering   \subfloat{\includegraphics[width=0.45\textwidth]{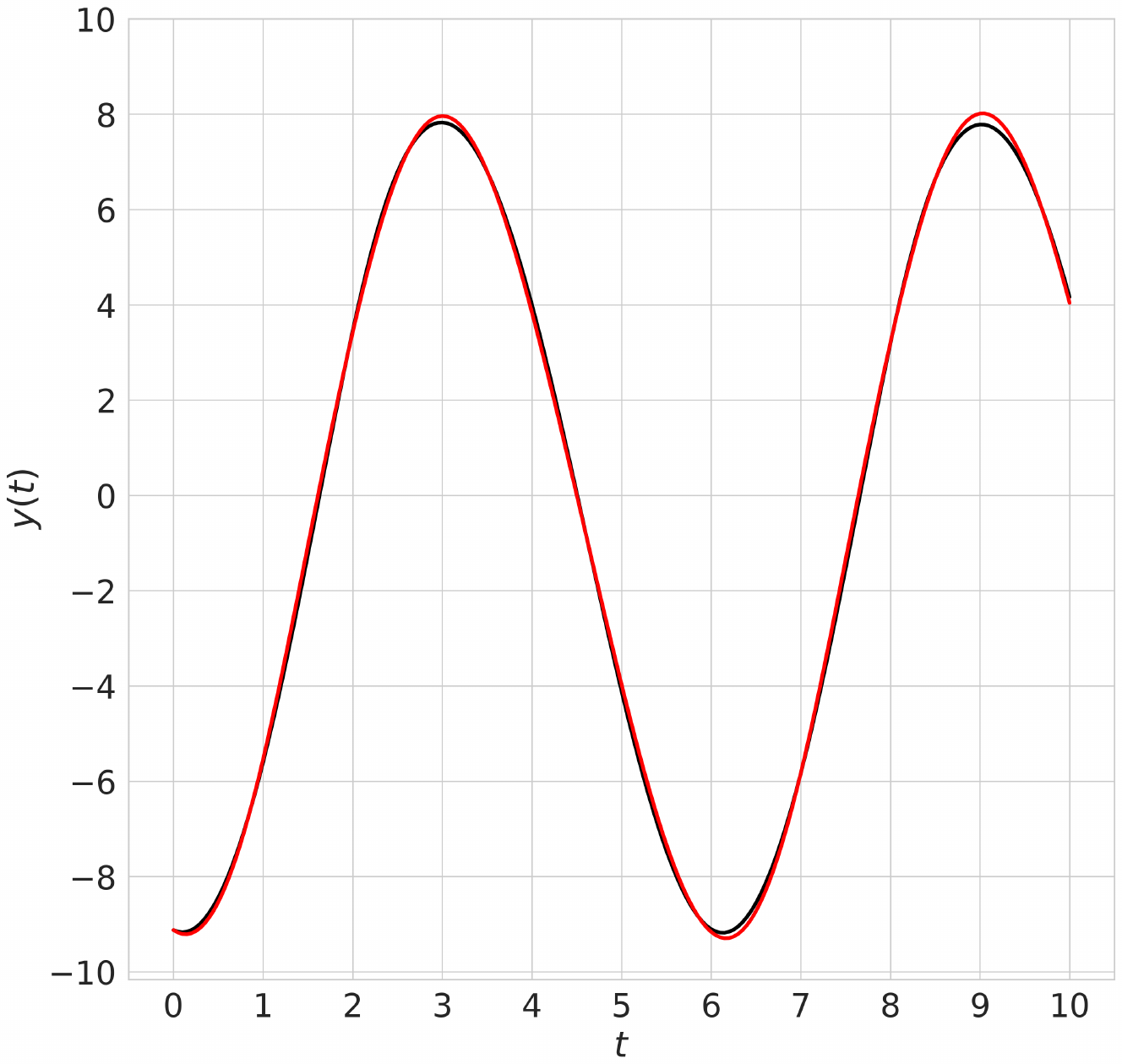}}
    \hspace{1mm}
     \subfloat{\includegraphics[width=0.45\textwidth]{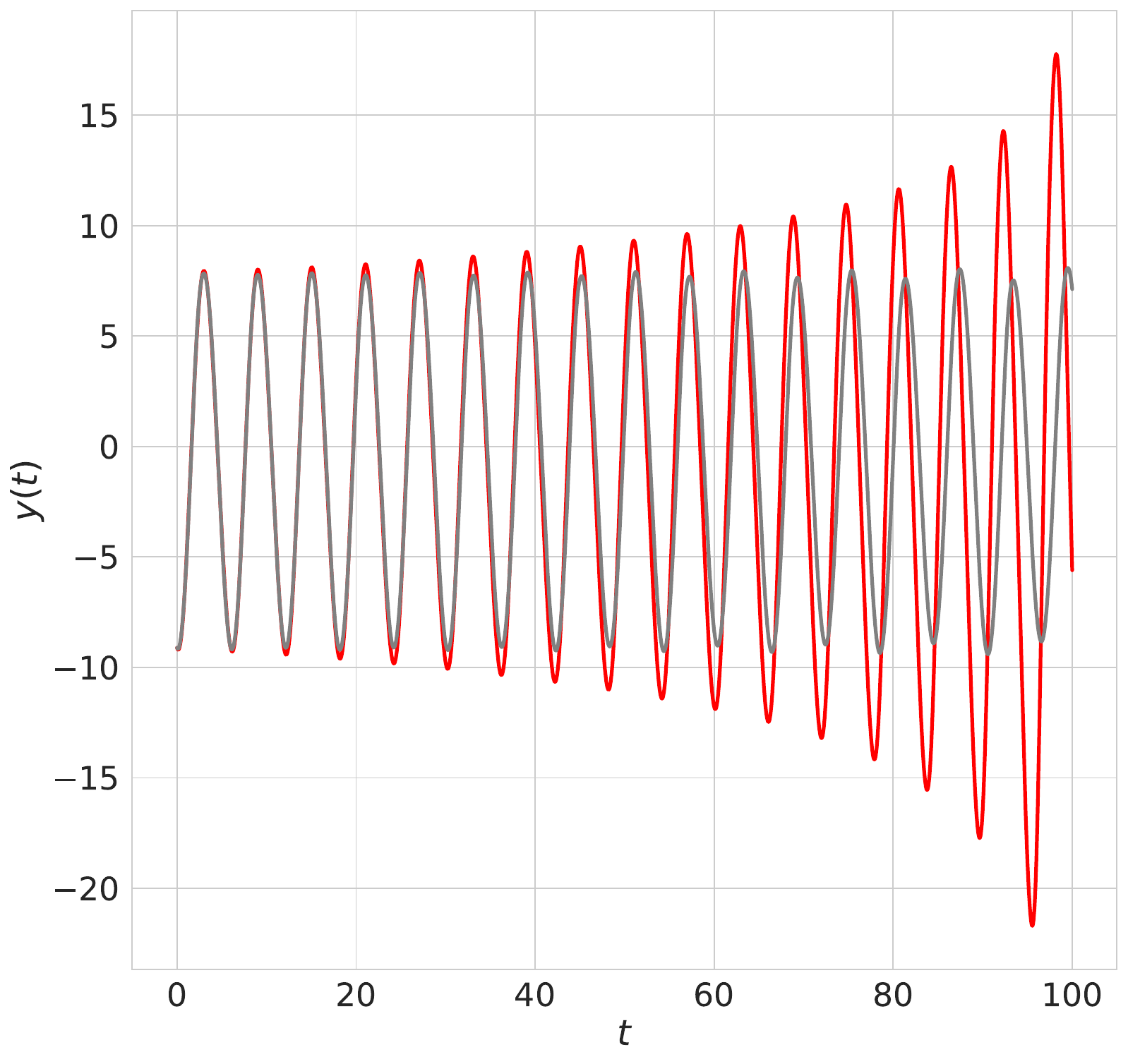}}
     \caption{\;Left: Training data in black along with the learned dynamics in red in the training time interval $[0, 10]$. Right: Target dynamics in gray along with the data generated by OS-net on the time interval $[0, 100]$.}\label{fig:ross_train_pred_noreg}
\end{figure}

We generate the training data by solving the R\"{o}ssler system using Matlab's ode45 solver with a time step of $0.001$ from $t = 0$ to $t = 10$. We then take snapshots of the simulated data every $50$ step and feed it to OS-net. We build OS-net using the Runge-Kutta 4 (RK4) algorithm with a time step of $0.005$. We chose the snake activation function $x + \frac{1}{0.2}\sin^2(0.2x)$ and set the number of nodes in the hidden layer to be $2 \times 16$. We set $\alpha = 0.07$ in Equation \myref{eq:loss_final} and use $10$ epochs.\\
Figure\myref{fig:ross_train_pred} (left) shows the training output for the $y$ component. OS-net was able to learn the dynamics accurately by the end of training. The norm of $J_a$ is approximately $0.99$ after training as recorded in Table\myref{tab:normJ}. In this case, Inequality\myref{eq:osnet_stab_cond} is not strictly enforced but the norm of the matrix $J_a$ is controlled enough so that OS-net renders stable orbits. The elements of the matrix $\Omega$ are concentrated in $[-0.7,\;0.7]$ as shown in Figure \ref{fig:colormap_ross}.\\
We validate OS-net by propagating a trajectory initialized $[ x_0, y_0, z_0] = [0, -9.1238 + 0.01 ,0]$ using the learned dynamics. Figure\myref{fig:ross_train_pred} (right) shows prediction using OS-net up to $t = 100$ and displays the accuracy of this prediction when compared to the correct dynamics.  We assess the stability of OS-net by propagating the trajectory to $t = 10000$. OS-net converges to a stable period-1 attractor while the R\"{o}ssler system converges to a period-2 one as showcased in Figure\myref{fig:ross_osnet}. 
\paragraph{Ablation study} We compare OS-net with a network obtained by keeping the same architecture and settings as in Section\myref{sec:ross6} but with the regularization in Equation\myref{eq:loss_final} switched off. The left side of Figure\myref{fig:ross_train_pred_noreg} shows that training was successful while the right side shows the dynamics learned by the unregularized network diverge from the true dynamics in the time interval $[0, 10]$. This shows the role of the regularization term in stabilizing the dynamics learned by OS-net.  \\

\subsubsection{c = 18, chaotic behavior}

We now set $c = 18$ and initialize the trajectory at $[ x_0, y_0, z_0] = [0, -22.9049 ,0]$. The R\"{o}ssler system displays a chaotic behavior in this regime. We generate the training data as before but take snapshots every $10$ steps. For OS-net, we use RK4 with a step size of $0.005$ and $x + \sin(x)$ as an activation function. The hidden layer size is $2 \times 32$ and the penalty coefficient $\alpha = 2$.\\ 
 \begin{figure}[htb]
    \centering   \subfloat{\includegraphics[width=0.45\textwidth]{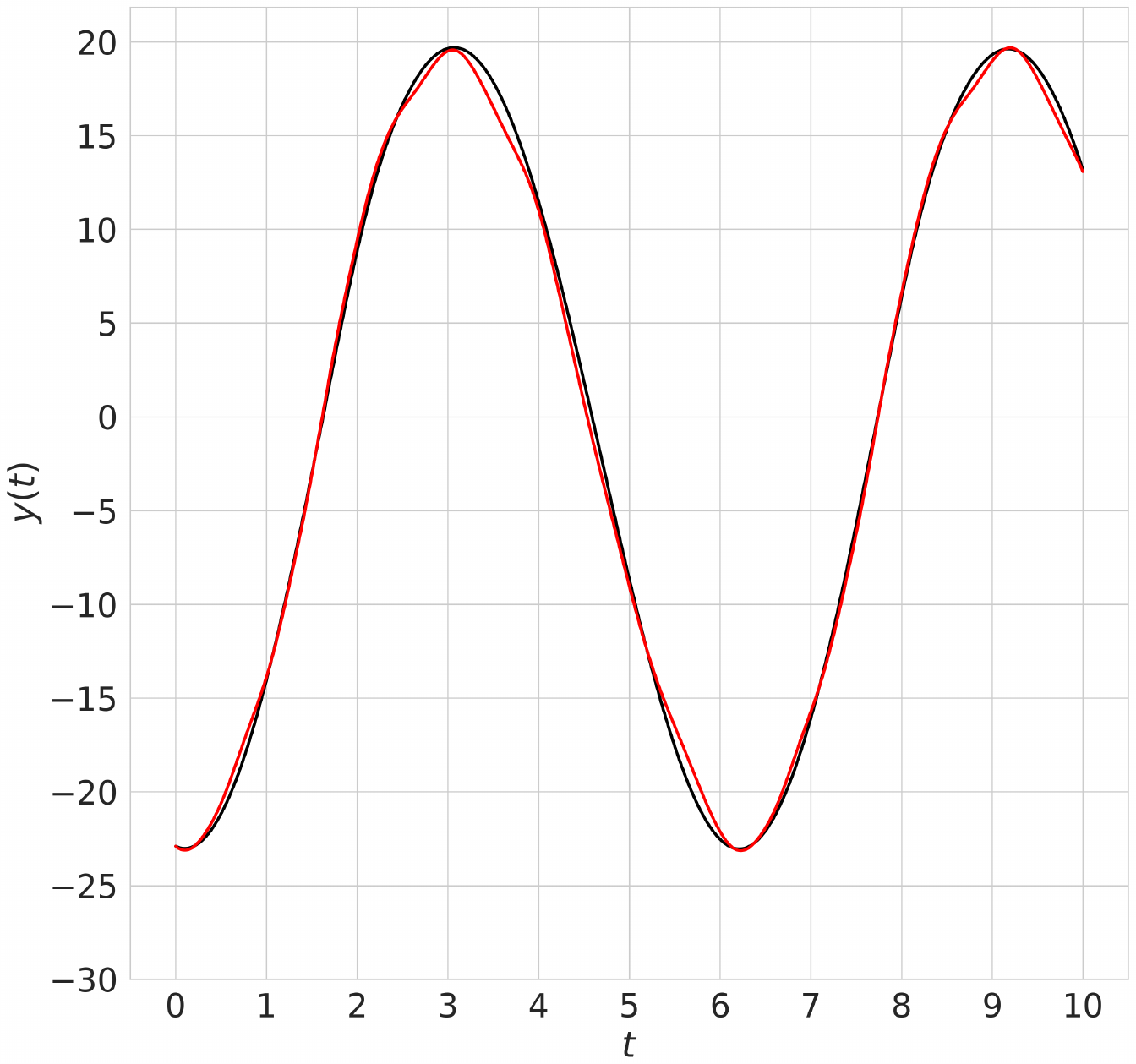}}
    \hspace{1mm}
     \subfloat{\includegraphics[width=0.45\textwidth]{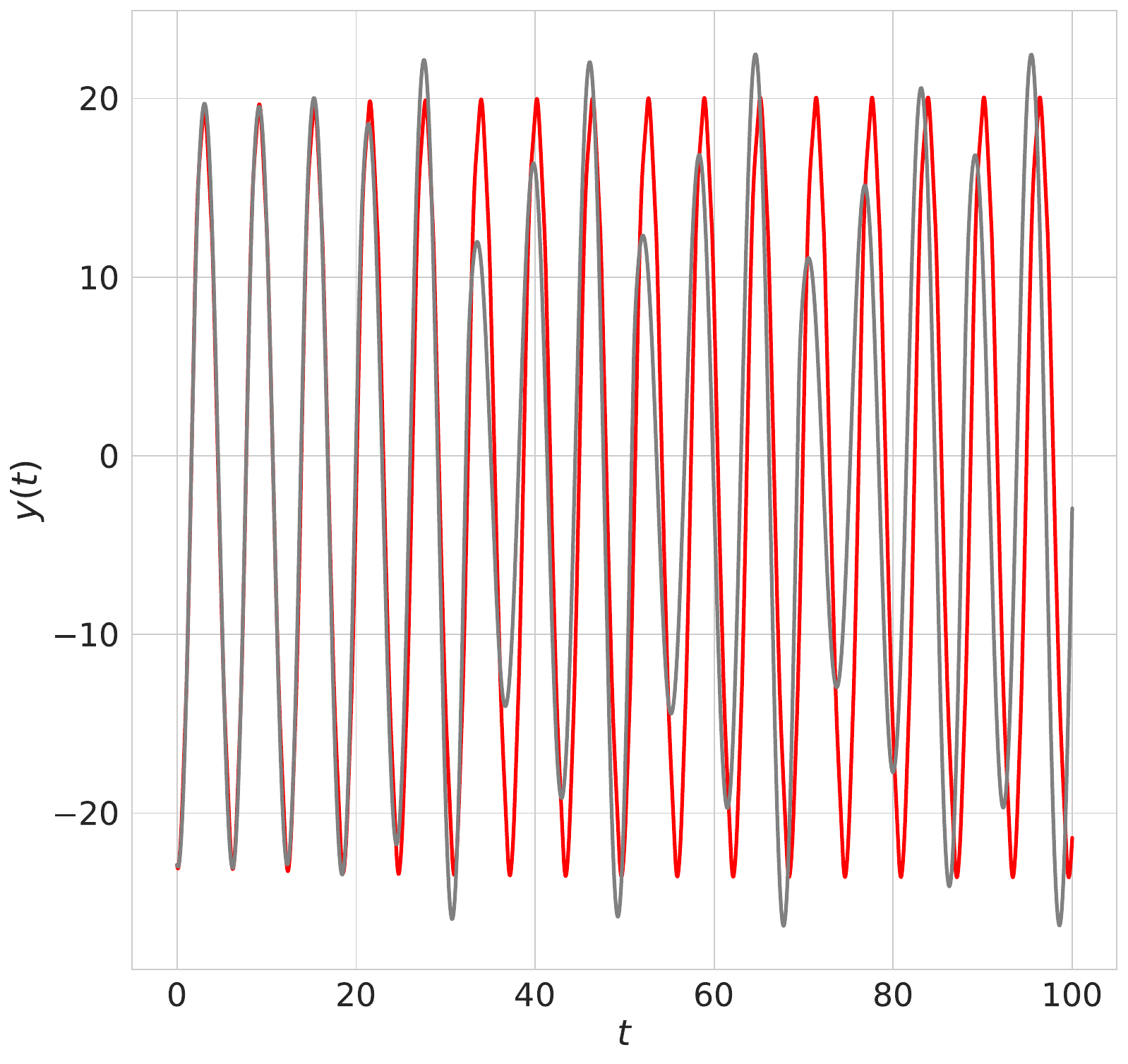}}
     \caption{\;Left: Training data in black along with the learned dynamics in red in the training time interval $[0, 10]$. Right: Target dynamics in gray along with the data generated by OS-net on the time interval $[0, 100]$.}\label{fig:ross_train_pred_18}
\end{figure}

 \begin{figure}[htb]
    \subfloat{\includegraphics[width=0.45\textwidth]{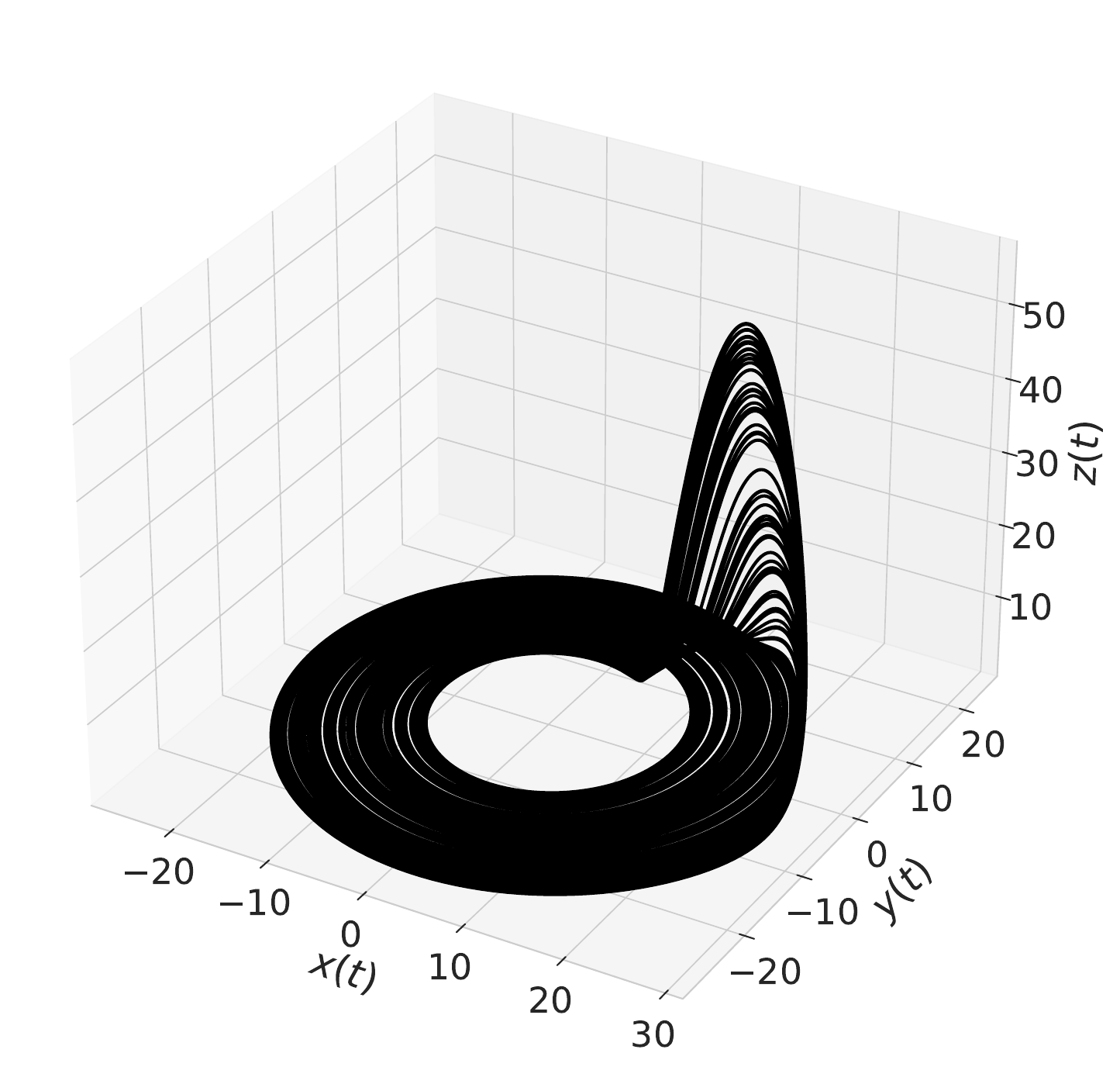}}
\hspace{2mm}
     \subfloat{\includegraphics[width=0.45\textwidth]{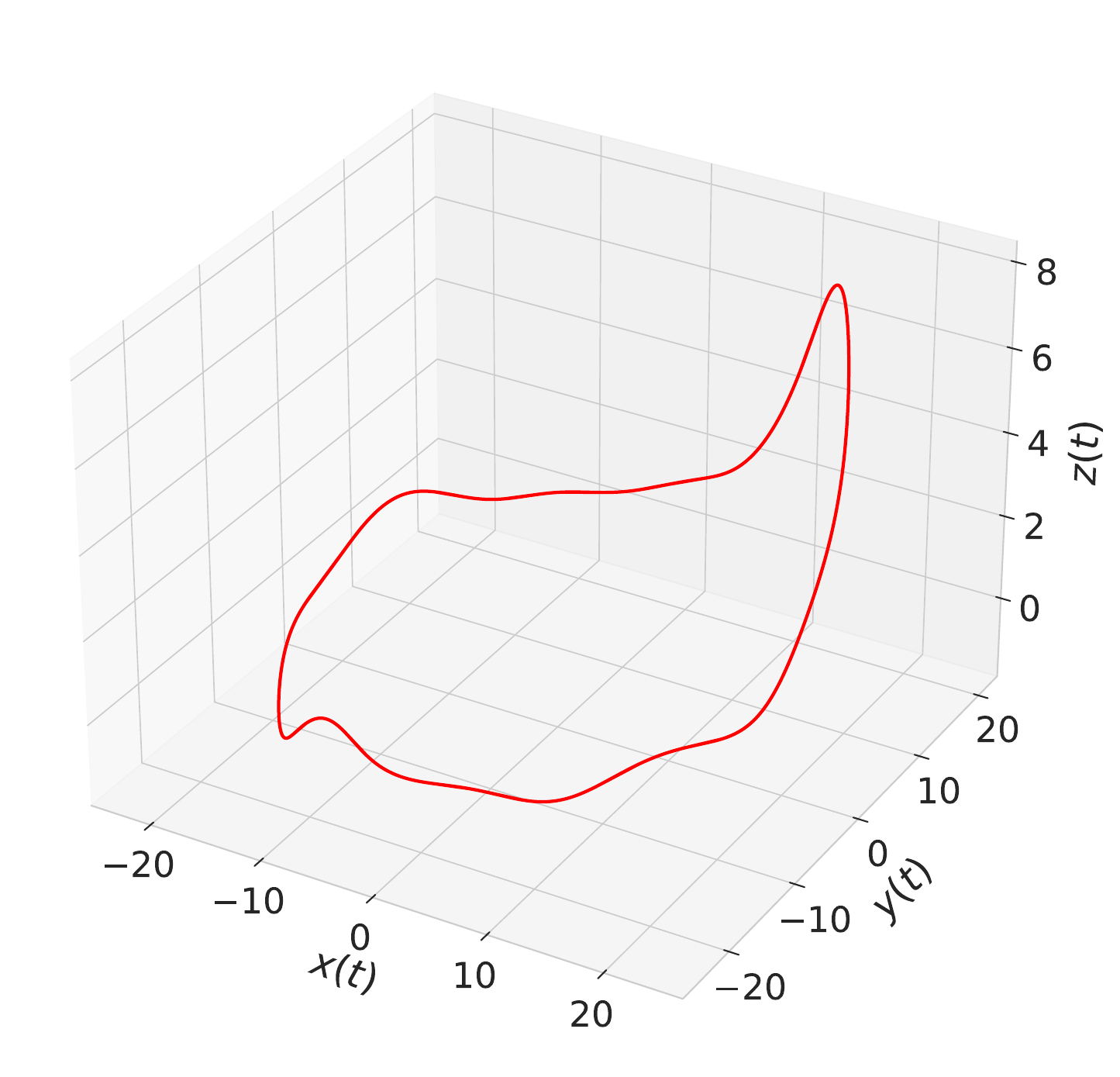}}
     \caption{\;Left: Chaotic R\"{o}ssler attractor. Right: Stable period-1  OS-net attractor }\label{fig:ross_c18}
    \centering
\end{figure}
Figure\myref{fig:ross_train_pred_18} shows the training output and confirms the ability of OS-net to learn the target dynamics. We then use the learned dynamics to generate a trajectory starting at $[ x_0, y_0, z_0] = [0, -22.9049 + 0.01 ,0]$. Since we are dealing with a chaotic system, the learned dynamics should not be expected to reproduce the training data \cite{Bramburger2020}. Figure\myref{fig:ross_train_pred_18} shows that OS-net was able to track the chaotic system up to $t \approx 30$. The norm of the matrix $J_a$ was approximately $0.6318$ at the end of training as recorded in Table\myref{tab:normJ}. Furthermore, the elements of the matrix $\Omega$ were concatenated between $-0.5$ and $+0.5$.  Figure\myref{fig:ross_c18} displays the chaotic R\"{o}ssler system and the stable attractor obtained by propagating OS-net's learned dynamics from $t = 0$ to $t = 10000$. 

\subsection{Simplest quadratic (Sprott's) chaotic flow}
\begin{figure}[htb]
    \subfloat{\includegraphics[width=0.45\textwidth]{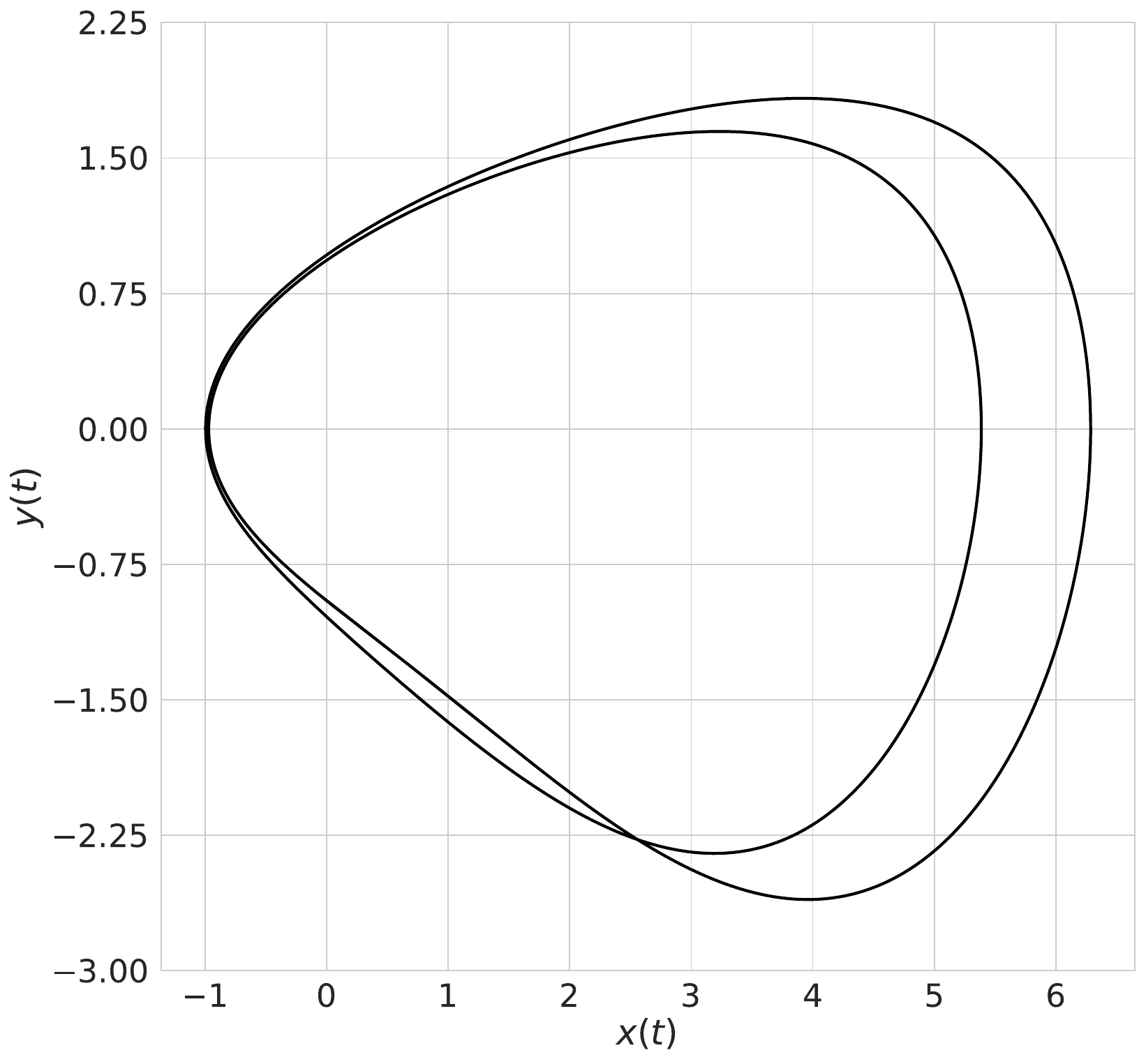}}
\hspace{2mm}
     \subfloat{\includegraphics[width=0.45\textwidth]{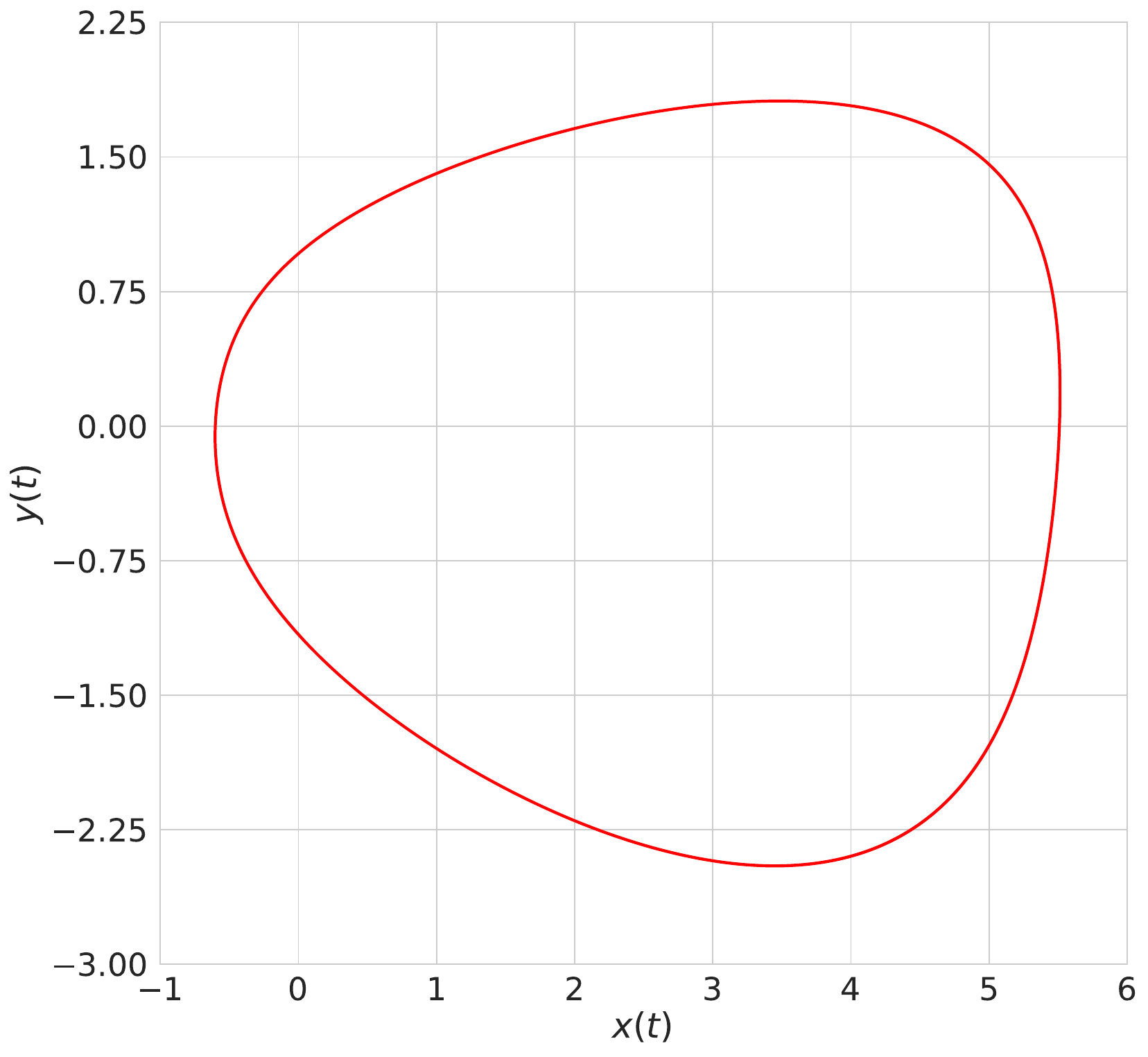}}
     \caption{\;Left: Sprott's period-2 attractor. Right: stable period-1 OS-net attractor}\label{fig:sprott_orbit}
    \centering
\end{figure}

We consider the following system 
\begin{align}\label{eq:sprott}
    \dot{x} &= y  \nonumber \\ 
    \dot{y} &= z \\
    \dot{z} &=-\nu z  -x + y^2 \nonumber
\end{align}
where $\nu \in \mathbf{R}$. This system was introduced in \cite{Sprott1997} and also has period doubling bifurcations as $\nu$ varies. Here we set $\nu = 2.1$ which yields a peiod-2 attractor for Equation\myref{eq:sprott}.\\
We initialize the trajectory at $[ x_0, y_0, z_0] =  [5.7043, 0.0 ,-2.12778]$ and solve the system using ode45 on the time interval  $[0, 15]$ with a step size of $0.001$. We then take snapshots every $10$ step and use the data for training. OS-net is solved using RK4 with a step size of $0.01$ and $x+\frac{1}{0.3}\sin^2(0.3 x)$ as an activation function. The hidden layer has $2 \times 16$ nodes and the penalty coefficient $\alpha = 1$.\\
 \begin{figure}[htb]
    \centering   \subfloat{\includegraphics[width=0.45\textwidth]{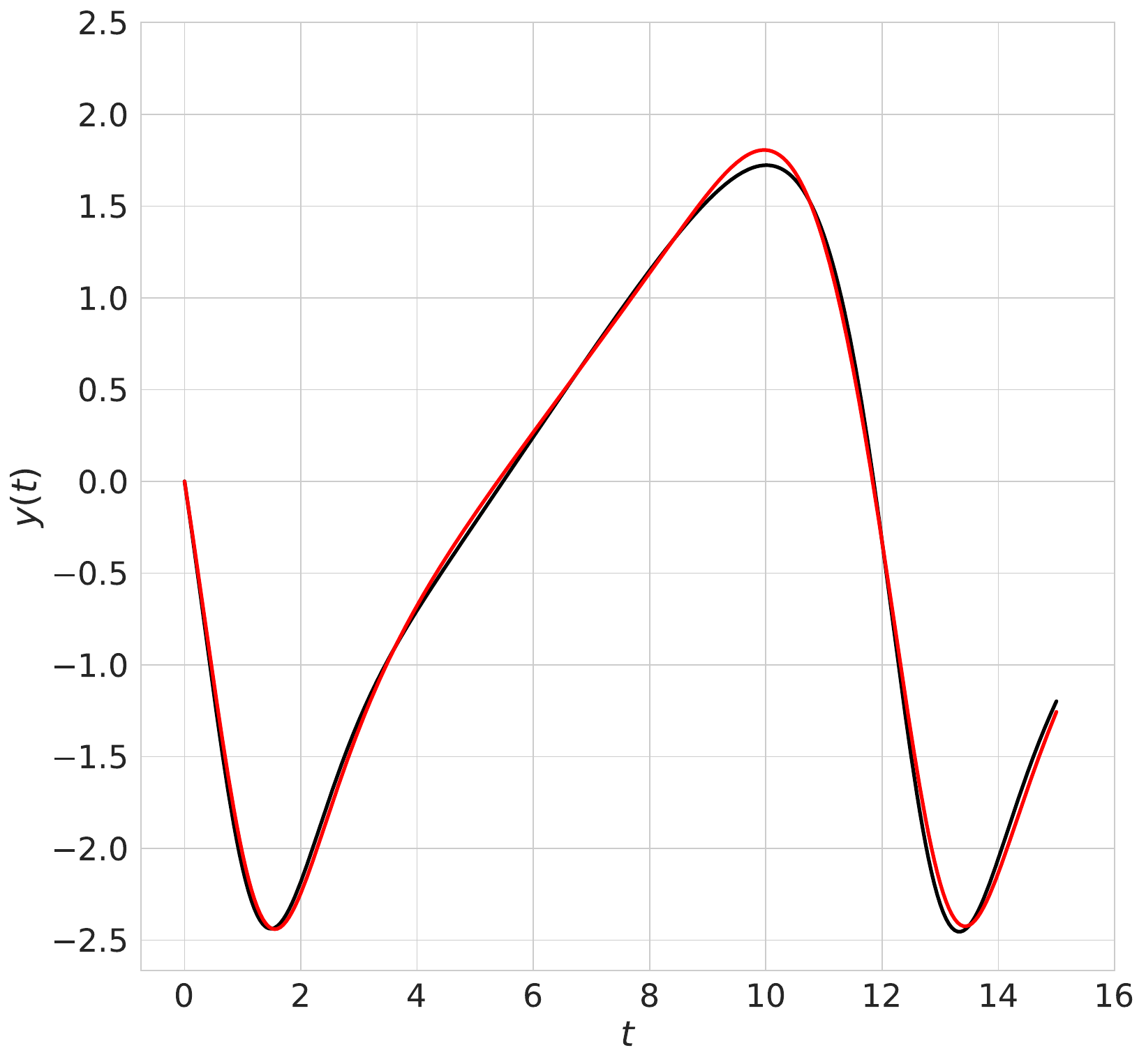}}
    \hspace{1mm}
     \subfloat{\includegraphics[width=0.45\textwidth]{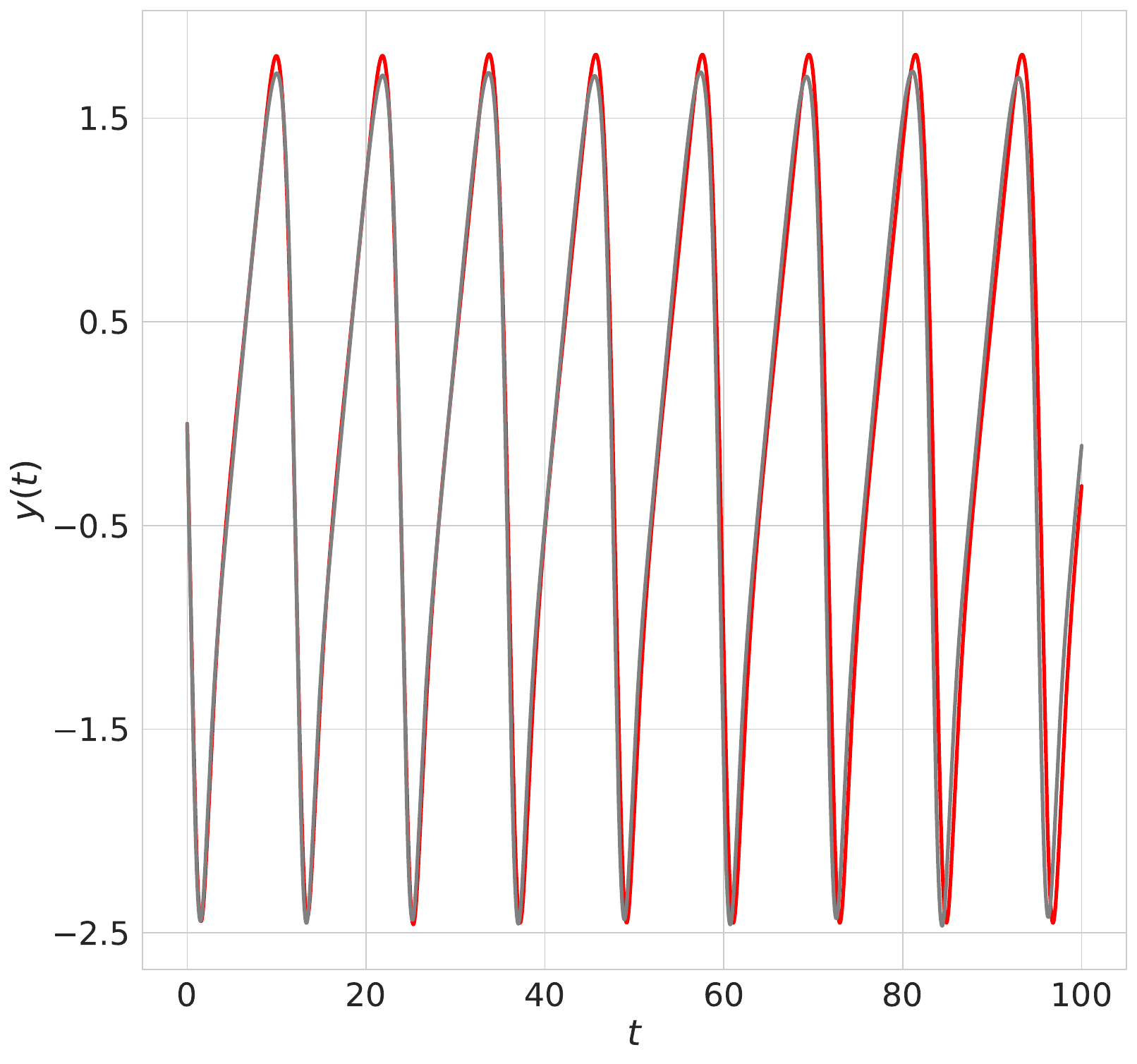}}
     \caption{\;Left: Training data in black along with the learned dynamics in red in the training time interval $[0, 15]$. Right: Target dynamics in gray along with the data generated by OS-net on the time interval $[0, 100]$.}\label{fig:sprott_train_pred}
\end{figure}
We show in Figure\myref{fig:sprott_train_pred} (left) the dynamics learned by OS-net for the $y$ component after $20$ epochs. Figure\myref{fig:sprott_train_pred} (right) also shows how well OS-net tracks the original system in the interval $t = 0$ to $t = 100$. In this case, the norm of the matrix $J_a$ was approximately $8e-3$ and the elements of the matrix $\Omega$ are in the interval $[-0.7,\;0.7]$. Inequality\myref{eq:osnet_stab_cond} is strictly enforced here. 
We then assess the stability of the learned dynamics by generating a trajectory starting at $[ x_0, y_0, z_0] =  [5.7043 + 0.01, 0.0,-2.12778]$ and evolving it from $t = 0$ to $t = 10000$. Figure\myref{fig:sprott_orbit} shows the period-2 attractor of the original system and the stable period-1 OS-net orbit. \\

\paragraph{Note}\label{sec:limitation}
The current implementation of OS-net uses the adjoint method presented in \cite{Chen2018} which accumulates numerical errors when integrating backward. We circumvent that by using RK4 with a small step size. This results in a computationally expensive implementation that can be improved using the methods proposed in \cite{Ott2021, zhuang2020, Zhang2022} that we plan on incorporating into OS-net in the future.

\section{Conclusion}\label{sec:conc}
We have presented a new family of stable neural network architectures for periodic dynamical data. The proposed architecture is a particular case of NODES with dynamics represented by a shallow neural network. We leveraged well-grounded ode theory to propose a new regularization scheme that controls the norm of the product of the weight matrices of the network. We have validated our theory by learning the R\"{o}ssler and Sprott's systems in different regimes including a chaotic one. In all the regimes considered, OS-net was able to track the exact dynamics and converge to a stable period-1 attractor. That indicates that OS-net is a promising network architecture that can handle highly complex dynamical systems.
In the future, we aim at controlling the parameters of the systems of interest by incorporating them into the state vectors that OS-net aims at learning. Additionally, we plan on using OS-net to learn and monitor the orbits of celestial objects that have short orbital periods such as certain exoplanets or three-body systems like Mars-Phobos. This extension of OS-net's applications holds great potential in providing a broader range of stable periodic orbits for the design of spatial missions.

\section*{Acknowledgments}
This work was supported by the U.S. Department of Energy, Office of Science, Office of Advanced Scientific Computing Research, Scientific Discovery through Advanced Computing (SciDAC) program and the FASTMath Institute under Contract No. DE-AC02-06CH11357 at Argonne National Laboratory. Government License. \\
{\bf Government License.}  The submitted manuscript has been created by
UChicago Argonne, LLC, Operator of Argonne National Laboratory
(``Argonne''). Argonne, a U.S. Department of Energy Office of Science
laboratory, is operated under Contract No. DE-AC02-06CH11357. The
U.S. Government retains for itself, and others acting on its behalf, a
paid-up nonexclusive, irrevocable worldwide license in said article to
reproduce, prepare derivative works, distribute copies to the public,
and perform publicly and display publicly, by or on behalf of the
Government.  The Department of Energy will provide public access to
these results of federally sponsored research in accordance with the
DOE Public Access
Plan. http://energy.gov/downloads/doe-public-access-plan.
\bibliographystyle{plain}
\bibliography{OSnet}

\begin{thebibliography}{10}

\bibitem{Bach2011}
F.~Bach, R.~Jenatton, J.~Mairal, and G.~Obozinski.
\newblock Optimization with sparsity-inducing penalties.
\newblock {\em Foundations and Trends in Machine Learning}, 4, 2011.

\bibitem{Bramburger2020}
J.~J. Bramburger and J.~N. Kutz.
\newblock Poincaré maps for multiscale physics discovery and nonlinear floquet
  theory.
\newblock {\em Physica D: Nonlinear Phenomena}, 408:132479, 2020.

\bibitem{Bramburger2021}
J.~J. Bramburger, J.~N. Kutz, and S.~L. Brunton.
\newblock Data-driven stabilization of periodic orbits.
\newblock {\em IEEE Access}, 9:43504--43521, 2021.

\bibitem{Brunton2016}
S.~L. Brunton, J.~L. Proctor, and N.~Kutz.
\newblock Discovering governing equations from data by sparse identification of
  nonlinear dynamical systems.
\newblock {\em Proceedings of the National Academy of Sciences of the United
  States of America}, 2016.

\bibitem{Chen2018}
T.~Q. Chen, Y.~Rubanova, J.~Bettencourt, and D.~Duvenaud.
\newblock Neural ordinary differential equations.
\newblock {\em 32nd Conference on Neural Information Processing Systems}, 2018.

\bibitem{Devaney2003}
R.~L. Devaney.
\newblock {\em An Introduction To Chaotic Dynamical Systems}.
\newblock CRC Press, 2003.

\bibitem{Dupont2019}
E.~Dupont, A.~Doucet, and Y.~W. Teh.
\newblock Augmented neural odes.
\newblock In H.~Wallach, H.~Larochelle, A.~Beygelzimer, F.~d\textquotesingle
  Alch\'{e}-Buc, E.~Fox, and R.~Garnett, editors, {\em Advances in Neural
  Information Processing Systems}, volume~32. Curran Associates, Inc., 2019.

\bibitem{ThetaGPU}
Argonne Leadership~Computing Facility.
\newblock Theta/thetagpu supercomputer.
\newblock \url{https://www.alcf.anl.gov/alcf-resources/theta}.

\bibitem{Haber2018}
E.~Haber and L.~Ruthotto.
\newblock Stable architectures for deep neural networks.
\newblock {\em Inverse Problems}, 34(1):014004, dec 2017.

\bibitem{horn2012matrix}
Roger~A Horn and Charles~R Johnson.
\newblock {\em Matrix analysis}.
\newblock Cambridge university press, 2012.

\bibitem{Krein1983}
M.~G. Krein.
\newblock Foundations of the theory of $\lambda$-zones of stability of a
  canonical system of linear differential equations with periodic coefficients.
\newblock {\em American Mathematical Society Translations}, 120, 1983.

\bibitem{KreinJakub1983}
M.~G. Krein and V.~A. Jakubovic.
\newblock Four papers on ordinary differential equations.
\newblock {\em American Mathematical Society Translations}, 120, 1983.

\bibitem{Ngom2021}
M.~Ngom and O.~Marin.
\newblock Fourier neural networks as function approximators and differential
  equation solvers.
\newblock {\em Statistical Analysis and Data Mining: The ASA Data Science
  Journal}, 14(6):647--661, 2021.

\bibitem{Ott2021}
K.~Ott, P.~Katiyar, P.~Hennig, and M.~Tiemann.
\newblock Resnet after all: Neural {\{}ode{\}}s and their numerical solution.
\newblock In {\em International Conference on Learning Representations}, 2021.

\bibitem{Parascandolo2017}
G.~Parascandolo, H.~Huttunen, and T.~Virtanen.
\newblock Taming the waves: sine as activation function in deep neural
  networks.
\newblock {\em https://openreview.net/pdf?id=Sks3zF9eg}, 2017.

\bibitem{Pontryagin1987}
L.~S. Pontryagin.
\newblock {\em Mathematical Theory of Optimal Processes}.
\newblock CRC Press, 1987.

\bibitem{raissi2019Pinn}
M.~Raissi, P.~Perdikaris, and G.~E. Karniadakis.
\newblock Physics-informed neural networks: A deep learning framework for
  solving forward and inverse problems involving nonlinear partial differential
  equations.
\newblock {\em Journal of Computational Physics}, 2019.

\bibitem{Roy1992}
R.~Roy, T.~W. Murphy, T.~D. Maier, Z.~Gills, and E.~R. Hunt.
\newblock Dynamical control of a chaotic laser: Experimental stabilization of a
  globally coupled system.
\newblock {\em Phys. Rev. Lett.}, 68:1259--1262, Mar 1992.

\bibitem{Ruthotto2019}
L.~Ruthotto and E.~Haber.
\newblock Deep neural networks motivated by partial differential equations.
\newblock {\em Journal of Mathematical Imaging and Vision volume}, 2019.

\bibitem{Rossler1976397}
O.~E. Rössler.
\newblock An equation for continuous chaos.
\newblock {\em Physics Letters A}, 57(5):397--398, 1976.

\bibitem{Schaeffer2017}
H.~Schaeffer, G.~Tran, and R.~Ward.
\newblock Learning dynamical systems and bifurcation via group sparsity, 2017.

\bibitem{Schiff1994}
S.~J. Schiff, D.~H. Jerger, K.and~Duong, M.~L. Chang, T.and~Spano, and W.~L.
  Ditto.
\newblock Controlling chaos in the brain.
\newblock {\em Nature}, 370:615--620, 1994.

\bibitem{Silvescu1999}
Adrian Silvescu.
\newblock Fourier neural networks.
\newblock {\em International joint conference on neural networks}, 1999.

\bibitem{Sirignano2018}
J.~Sirignano and K.~Spiliopoulos.
\newblock Dgm: A deep learning algorithm for solving partial differential
  equations.
\newblock {\em Journal of Computational Physics}, 2018.

\bibitem{Sprott1997}
J.~C. Sprott.
\newblock Simplest dissioative chaotic flow.
\newblock {\em Physics Letters A}, 228, 1997.

\bibitem{Sprott2003}
J.~C. Sprott.
\newblock {\em Chaos and time-series analysis}.
\newblock Oxford University Press, 2003.

\bibitem{Sun2020}
Y.~Sun and H.~Zhang, L.and~Schaeffer.
\newblock A memory-efficient neural ordinary differential equation framework
  based on high-level adjoint differentiation.
\newblock {\em Proceedings of Machine Learning Research}, 107:352--372, 2020.

\bibitem{msurtsukov_github}
M.~Surtsukov.
\newblock Pytorch implementation of neural ordinary differential equations.
\newblock \url{https://github.com/msurtsukov/neural-ode}, 2019.

\bibitem{Teschl2012}
G.~Teschl.
\newblock {\em Ordinary Differential Equations and Dynamical Systems}.
\newblock AMS, 2012.

\bibitem{Wiesel1993}
W.~Wiesel and W.~Shelton.
\newblock Modal control of an unstable periodic orbit.
\newblock {\em Journal of the Astronautical Sciences}, 31:63--76, 1983.

\bibitem{Wolfe1969}
P.~Wolfe.
\newblock Convergence conditions for ascent methods.
\newblock {\em SIAM Review}, 11(2):226--235, 1969.

\bibitem{Yakubovich1975}
V.~A. Yakubovich and V.~M. Starzhinskii.
\newblock {\em Linear Differential Equations with Periodic Coefficients}.
\newblock Halsted Press, 1975.

\bibitem{Yan2020}
H.~Yan, J.~Du, V.~Tan, and J.~Feng.
\newblock On robustness of neural ordinary differential equations.
\newblock In {\em International Conference on Learning Representations}, 2020.

\bibitem{Zhang2022}
H.~Zhang and W.~Zhao.
\newblock A memory-efficient neural ordinary differential equation framework
  based on high-level adjoint differentiation.
\newblock {\em IEEE Transactions on Artificial Intelligence}, pages 1--11,
  2022.

\bibitem{zhuang2020}
J.~Zhuang, X.~Dvornek, N.and~Li, S.~Tatikonda, X.~Papademetris, and J.~Duncan.
\newblock Adaptive checkpoint adjoint method for gradient estimation in neural
  {ODE}.
\newblock In Hal~Daumé III and Aarti Singh, editors, {\em Proceedings of the
  37th International Conference on Machine Learning}, volume 119 of {\em
  Proceedings of Machine Learning Research}, pages 11639--11649. PMLR, 13--18
  Jul 2020.

\bibitem{Zhumekonov2019}
A.~Zhumekenov, M.~Uteuliyeva, O.~Kabdolov, R.~Takhanov, Z.~Assylbekov, and
  A.~J. Castro.
\newblock Fourier neural networks: {A} comparative study.
\newblock {\em http://arxiv.org/abs/1902.03011}, 2019.

\bibitem{Ziyin2020}
L.~Ziyin, T.~Hartwig, and M.~Ueda.
\newblock Neural networks fail to learn periodic functions and how to fix it.
\newblock In {\em Advances in Neural Information Processing Systems},
  volume~33, 2020.

\end{thebibliography}

\appendix
\section{Orbits of dynamical systems}\label{appendixA}
In this section, we recall fundamental results for the stability of periodic orbits of equations of the form
\begin{equation}\label{eq:simple_ode}
    \dot{\bm x} = f(x),\;\; \bm x(0) = \bm x_0.
\end{equation}
where $f \in \mathrm{C}^k(\mathrm{M}, \mathrm{R}^n)$, $k \geq 1$ and $\mathcal{M}$ an open subset of $\mathrm{R}^n$.

Chaos theory and the stability of periodic orbits has been widely studied in the literature \cite[Chapter~12]{Teschl2012}, \cite{Sprott2003}. They are important concepts in fields like celestial mechanics, biology and chemistry. 
An important tool to study periodic orbits is the Poincar\'e or return map. Let $\phi(t, \bm x_0)$ be a periodic solution of Equation \myref{eq:simple_ode} with period $T$. Let 
$I_x$ be the maximal interval where $\phi$ is defined. Denoting $$\mathrm{W} = \cup_{x \in \mathcal{M}} I_x \times \{x\} \subseteq	\mathrm{R} \times \mathrm{M},$$ the flow of Equation\myref{eq:simple_ode} is defined to be the map
$$\Phi:\;\mathrm{W} \rightarrow \mathrm{M},\;\; (t,x) \mapsto \phi(t,x).$$ Let $\gamma(\bm x_0)$ be the associated periodic orbit. The Poincar\'e map is defined as 
\begin{equation}\label{eq:poincare}
    P_{\Sigma}(\bm y) = \Phi(\tau(\bm y),\bm y)
\end{equation}
where $\Sigma$ is a transversal submanifold of codimension one containing one value $\bm x_0$ from the periodic orbit $\gamma(\bm x_0)$, and $\tau \in \mathcal{C}^k(U)$ such that $\tau(x_0) = T$ and $U$ a neighborhood of $x_0$ such that $\forall\; y \in U$, $\Phi(\tau(\bm y),\bm y) \in \Sigma$. It has been proven that the stability of periodic orbits is directly connected to the stability of $\bm x_0$ as a fixed point of the return map $P_{\Sigma}$. More precisely, we have the following theorem \cite[Chapter 12]{Teschl2012}
\begin{theorem2}
Suppose $f \in \mathrm{C}^k$ has a periodic orbit $\gamma(\bm x_0)$. If all eigenvalues of the derivative of the Poincare map $DP_{\Sigma}$ at $\bm x_0$ lie inside the unit circle then the periodic orbit is asymptotically stable.
\end{theorem2}

It is however generally difficult to compute Poincar\'e maps and their derivatives explicitly. Fortunately, it was proven in \cite{Teschl2012} that the eigenvalues of the derivative of  Poincare map $DP_{\Sigma}$ at $\bm x_0$ plus the single value $1$ coincide with the eigenvalues of the monodromy matrix (see Definition\myref{def:monodrony}) of the first variational (FV) equation associated with Equation \myref{eq:simple_ode}
\begin{equation}\label{eq:fvar_app}
    \Dot{\bm y} = A(t)\bm y, \;\; \bm y(t_0) = \bm I,\;\;\bm A(t) = d (f(\bm x))_{(\Phi(t,\bm x_0))}\;\; \text{and } \bm A(t+T) = \bm A(t).
\end{equation}
Therefore, by evaluating the stability of the FV equation associated with a periodic orbit one can assess its stability properties.
\section{Stability of Linear Canonical systems}\label{appendixB}
In addition to the definitions given in the main body of this article, we give more definitions and results that would allow us to prove Theorem \ref{theo:krein}.
\begin{definition}\label{def:monodrony}
    The monodromy $\bm U(T)$ of a periodic linear system $\dot(x) = \bm A(t)x,\;\; \bm A(t+T) = \bm A(t)$ is $$\bm U(T) = \bm \Pi(T,t_0)$$ where $\bm \Pi(t,t_0)$ is the principal matrix solution of the system i.e. $\bm \Pi(t,t_0)$  solves the initial value problem $$\dot{\bm \Pi}(t,t_0) =\bm A(t) \bm \Pi(t,t_0),\;\; \bm \Pi(t,t_0) = \bm I.$$ 
\end{definition}
\begin{definition}
    A matrix $\bm S$ is said to be J-unitary if $\bm U^* \bm J \bm U = \bm J$. In particular, the monodromy matrix $\bm U(T, \lambda)$ of Equation\myref{eq:krein_cano} is J-unitary.
\end{definition}
\begin{definition}
    Let $\bm H \in \mathrm{P}_n(T)$ and consider the boundary value problem (BVP) 
    \begin{align*}
       \dot{\bm y} = \lambda \bm J_m \bm H(t) \bm y,\;\; \bm y(T) = \bm \Xi \bm y(0)  
    \end{align*}
where $\bm \Xi$ is a J-unitary matrix. The characteristics values of this BVP are the roots (for $\lambda$) of the equation $$ det(\bm U(T, \lambda) - \bm \Xi) = 0.$$ 
\end{definition}

We now recall the following results from \cite[Theorem 6.1 and Theorem 6.2]{Krein1983}
\begin{theorem}
    If $\bm H \in \mathrm{P}_n(T)$, then the BVP
    \begin{align}
       \dot{\bm y} = \lambda \bm J_m \bm H(t) \bm y,\;\; \bm y(T) = -\bm y(0)  
    \end{align}
    has at least one positive and one negative characteristic value. Furthermore, if we denote $\Lambda_+$ the smallest positive characteristic value and $\Lambda_-$ the largest negative one, then the open interval $(\Lambda_-, \Lambda_+)$ belongs to the central zone of stability of Equation\myref{eq:krein_cano}.
\end{theorem}\label{theo:exis}
We recall the following theorem
\begin{theorem}[Krein]
    A real $\lambda$ belongs to the central zone of stability of an Equation \myref{eq:krein_cano} of positive type, if 
$$|\lambda| < 2 \mathcal{M}^{-1} (\bm C)$$
where $\bm C = \bm J_{m_{a}}\int_0^T \bm H_a(t)$. 
If $K$ is a matrix, $K_a$ is the matrix obtained by replacing the elements of $K$ by their absolute values.
\end{theorem}
\begin{proof}
Let $\bm x^+(t) = (x^+_1, \cdots,  x^+_n)$ be a non trivial solution of Equation\myref{eq:osnet_fvar} for $\Lambda_+$ such that $\bm x^+(t+T) = -\bm x^+(T)$. The existence of $\Lambda_+$ and $\bm x^+$ is assured by Theorem\myref{theo:exis}.

The corresponding system is $$ \dot{\bm x^+} = \Lambda_+ \bm J \bm H(t)  \bm x^+$$ where $\bm x^+(t+T) = -\bm x^+(t)$ and we set $v_j = \max_{0\leq t \leq T} |x^+_j(t)| = |x_j^+(\tau_j)|,\; j = 1\cdots 2m.$

Let $\bm A = \bm J \bm H(t)$,  we have
$$\dot{x^+}_j = \Lambda^+ \sum_k a_{jk}(t)x^+_k, \;\; j = 1\cdots 2m.$$
We integrate these equations from $\tau_j$ to $\tau_j +T$ to obtain
$$-2 x_j^+(\tau_j) =  \Lambda_+ \sum_k \int_{\tau_j}^{\tau_j + T}c_{jk}(t)x^+_k(t) dt$$
We now take the modulus and obtain
$$ 2v_j \leq \Lambda_+ \sum_k v_k \int_{\tau_j}^{\tau_j + T} |a_{jk}(t)| dt =\Lambda_+ \sum_k c_{jk}v_k$$
Where $ \bm C = (c_{ij})_{j=1,\cdots, n} = \bm J_a \int_{0}^{T} \bm H_a(t) dt $ and obtain
$$ \bm v \leq \frac{1}{2} \left(  \Lambda_+ J_a \int_{0}^{T} \bm H_a(t)dt \right) \bm v$$ where $\bm v = (v_1,\cdots, v_n)$. 
We now use the following lemma proven in \cite{Krein1983}:
\begin{lemma}
If for a nonzero matrix $\bm A = (a_{ij})_{i,j = 1,\cdots, n} $ with nonnegative elements there exists a nonzero vector $\bm v = (v_1, \cdots, v_n)$ with nonnegative coordinates such that $\bm v \leq \bm A \bm v,$ then $\mathcal{M}(\bm A) \geq 1$.
\end{lemma}
to state
$$\Lambda_+\frac{\mathcal{M}\left(\bm J_a \int_{0}^{T} \bm H_a(t)dt\right) }{2}\geq 1 \;\;\;\; \text{i.e.} \;\;\;\; \Lambda_+ \geq 2 \mathcal{M}^{-1}\left( \bm J_a \int_{0}^{T} \bm H_a(t)dt \right).$$
In a similar fashion, we obtain $- \Lambda_- \geq 2 \mathcal{M}^{-1}\left( \bm J_a \int_{0}^{T} \bm H_a(t)dt\right).$
Consequently,  $\lambda$ is in the central zone of stability if $$|\lambda| < 2 \mathcal{M}^{-1}\left(\bm J_a \int_{0}^{T} \bm H_a(t)dt\right).$$
Hence, if $1 < 2 \mathcal{M}^{-1}\left(\bm J_a \int_{0}^{T} \bm H_a(t)dt\right)$, then Equation\myref{eq:krein_cano} is stable.
\end{proof}

\end{document}